\documentclass{article}

\usepackage{amssymb,latexsym,amsmath}

\usepackage{graphicx}

\begin{document}

\newcommand{\bfi}{\bfseries\itshape}

\makeatletter

\@addtoreset{figure}{section}

\def\thefigure{\thesection.\@arabic\c@figure}

\def\fps@figure{h, t}

\@addtoreset{table}{bsection}

\def\thetable{\thesection.\@arabic\c@table}

\def\fps@table{h, t}

\@addtoreset{equation}{section}

\def\theequation{\thesubsection.\arabic{equation}}

\makeatother

\newtheorem{theorem}{Theorem}[section]

\newtheorem{proposition}[theorem]{Proposition}

\newtheorem{lema}[theorem]{Lemma}

\newtheorem{corollary}[theorem]{Corollary}

\newtheorem{definition}[theorem]{Definition}

\newtheorem{remark}[theorem]{Remark}

\newtheorem{exempl}{Example}[section]

\newenvironment{exemplu}{\begin{exempl}  \em}{\hfill $\square$

\end{exempl}}

\newcommand{\comment}[1]{\par\noindent{\raggedright\texttt{#1}

\par\marginpar{\textsc{Comment}}}}

\newcommand{\todo}[1]{\vspace{5 mm}\par \noindent \marginpar{\textsc{ToDo}}\framebox{\begin{minipage}[c]{0.95 \textwidth}

\tt #1 \end{minipage}}\vspace{5 mm}\par}

\newcommand{\ea}{\mbox{{\bf a}}}

\newcommand{\eu}{\mbox{{\bf u}}}

\newcommand{\ueu}{\underline{\eu}}

\newcommand{\ueo}{\overline{u}}

\newcommand{\oeu}{\overline{\eu}}

\newcommand{\ew}{\mbox{{\bf w}}}

\newcommand{\ef}{\mbox{{\bf f}}}

\newcommand{\eF}{\mbox{{\bf F}}}

\newcommand{\eC}{\mbox{{\bf C}}}

\newcommand{\en}{\mbox{{\bf n}}}

\newcommand{\eT}{\mbox{{\bf T}}}

\newcommand{\eL}{\mbox{{\bf L}}}

\newcommand{\eR}{\mbox{{\bf R}}}

\newcommand{\eV}{\mbox{{\bf V}}}

\newcommand{\eU}{\mbox{{\bf U}}}

\newcommand{\ev}{\mbox{{\bf v}}}

\newcommand{\eve}{\mbox{{\bf e}}}

\newcommand{\uev}{\underline{\ev}}

\newcommand{\eY}{\mbox{{\bf Y}}}

\newcommand{\eK}{\mbox{{\bf K}}}

\newcommand{\eP}{\mbox{{\bf P}}}

\newcommand{\eS}{\mbox{{\bf S}}}

\newcommand{\eJ}{\mbox{{\bf J}}}

\newcommand{\eB}{\mbox{{\bf B}}}

\newcommand{\eH}{\mbox{{\bf H}}}

\newcommand{\leb}{\mathcal{ L}^{n}}

\newcommand{\eI}{\mathcal{ I}}

\newcommand{\eE}{\mathcal{ E}}

\newcommand{\hen}{\mathcal{H}^{n-1}}

\newcommand{\eBV}{\mbox{{\bf BV}}}

\newcommand{\eA}{\mbox{{\bf A}}}

\newcommand{\eSBV}{\mbox{{\bf SBV}}}

\newcommand{\eBD}{\mbox{{\bf BD}}}

\newcommand{\eSBD}{\mbox{{\bf SBD}}}

\newcommand{\ecs}{\mbox{{\bf X}}}

\newcommand{\eg}{\mbox{{\bf g}}}

\newcommand{\paromega}{\partial \Omega}

\newcommand{\gau}{\Gamma_{u}}

\newcommand{\gaf}{\Gamma_{f}}

\newcommand{\sig}{{\bf \sigma}}

\newcommand{\gac}{\Gamma_{\mbox{{\bf c}}}}

\newcommand{\deu}{\dot{\eu}}

\newcommand{\dueu}{\underline{\deu}}

\newcommand{\dev}{\dot{\ev}}

\newcommand{\duev}{\underline{\dev}}

\newcommand{\weak}{\stackrel{w}{\approx}}

\newcommand{\mild}{\stackrel{m}{\approx}}

\newcommand{\lrightarrow}{\stackrel{L}{\rightarrow}}

\newcommand{\rrightarrow}{\stackrel{R}{\rightarrow}}

\newcommand{\strong}{\stackrel{s}{\approx}}

\newcommand{\weakdown}{\rightharpoondown}

\newcommand{\opg}{\stackrel{\mathfrak{g}}{\cdot}}

\newcommand{\opunu}{\stackrel{1}{\cdot}}
\newcommand{\opdoi}{\stackrel{2}{\cdot}}

\newcommand{\opn}{\stackrel{\mathfrak{n}}{\cdot}}
\newcommand{\opx}{\stackrel{x}{\cdot}}

\newcommand{\tr}{\ \mbox{tr}}

\newcommand{\Ad}{\ \mbox{Ad}}

\newcommand{\ad}{\ \mbox{ad}}

\renewcommand{\contentsname}{ }

\title{What is a space? Computations in emergent algebras and the front end visual system}

\author{Marius Buliga \\
\\
Institute of Mathematics, Romanian Academy, 
P.O. BOX 1-764, \\
 RO 014700, Bucure\c sti, Romania\\
{\footnotesize Marius.Buliga@imar.ro}}

\date{This version:  25.09.2010}

\maketitle

\begin{abstract}
With the help of link diagrams with decorated crossings, I explain computations in 
emergent algebras, introduced in \cite{buligairq}, as the  kind of computations done in the 
front end visual system. 
\end{abstract}

\vspace{.5cm}





\section{Computations in the front end visual system as a paradigm}

In mathematics "spaces" come in many flavours. There are vector spaces, affine spaces, symmetric spaces, groups and so on. We usually take such objects as the stage where the plot of reasoning is laid. 
But in fact what we use, in many instances, are properties of particular spaces which, I claim, can be seen as coming from a particular class of computations.  

There is though a "space" which is "given" almost beyond doubt, namely the physical space where we all live. But as it regards perception of this space, we know now that things are not so simple. 
As I am writing these notes, here in Baixo Gavea, my eyes are attracted by the wonderful complexity 
of a tree near my window. The nature of the tree is foreign to me, as are the other smaller beings growing on or around the tree. I can make some educated guesses about what they are: some are orchids, there is a smaller, iterated version of the big tree. However, somewhere in my brain,  at a very fundamental level, the visible space is constructed in my head, before the stage where I am capable of recognizing and naming the objects or beings that I see. I cite from Koenderink  \cite{koen}, p. 126: 
\vspace{.5cm}

"The brain can organize {\em itself} through information obtained via interactions with the physical world into 
an embodiment of geometry, it becomes a veritable {\em geometry engine}. [...] 

There may be a point in holding that many of the better-known brain processes are most easily understood in terms of differential geometrical calculations running on massively parallel processor arrays whose nodes can be understood quite directly in terms of multilinear operators (vectors, tensors, etc). In this view brain processes in fact are space."

\vspace{.5cm}
 
In the paper \cite{koen2} Koenderink, Kappers and van Doorn study the "front end visual system" starting from general invariance  principles.  They define the "front end visual system" as "the interface between the light field and those parts of the brain nearest to the transduction stage".  After explaining that 
the "exact limits of the interface are essentially arbitrary", the authors propose the following 
characterization of what the front end visual system does, section 1 \cite{koen2}: 
\begin{enumerate}
\item[1.] "the front end is a "machine" in the sense of a syntactical transformer; 
\item[2.] there is no semantics. The front end processes structure;
\item[3.] the front end is precategorical, thus -- in a way -- the front end does not compute anything; 
\item[4.] the front end operates in a bottom-up fashion. Top down commands based upon semantical interpretations are not considered to be part of the front end proper;
\item[5.] the front end is a deterministic machine; ... all output depends causally on the (total) input from the immediate past." 
\item[6.] "What is not explicitly encoded by the front end is irretrievably lost. Thus the front end should be 
universal (undedicated) and yet should provide explicit data structures (in order to sustain fast processing past the front end) without sacrificing completeness (everything of potential importance to the survival of the agent has to be represented somehow)."
\end{enumerate}
 
 The authors continue by explaining that in the brain there is an embodiment of a fiber bundle, with 
 the visual field  as the basis and with cortical hypercolumns as  fibers. There is therefore 
 a piece of hardware in the brain which allows every point of the visual field to carry a copy of the 
 tangent space. By this, they mean that the hardware of the brain, if fed with the image of an edge at a point in the visual field,   it is capable of representing it as a tangent vector at that point. 
 
 This is one of the basic things that the front end does. In general, the front end does things in a 
 massively parallel manner, thus it has to do it by working with local representations. I cite from the last paragraph of section 1, just before the beginning of sections 1.1 \cite{koen2}.

\vspace{.5cm}

"In a local representation one can do without extensive (that is spatial, or geometrical) properties and represent everything in terms of intensive properties. This obviates the need for explicit geometrical expertise. The local representation of geometry is the typical tool of differential geometry. ... The columnar organization of representation in primate visual cortex suggests exactly such a structure."

\vspace{.5cm}

So, the front end does perform a kind of a computation, although a  very particular one. It is not, 
at first sight, a logical, boolean type of computation. I think this is what Koenderink and coauthors want to say in  point 3. of  the characterization of the front end visual system. Nevertheless, if you want 
your visual system to be hacked such that to perform logical computations, see this \cite{changizi}.

 We may imagine that there is an abstract mathematical "front end" which, if fed with the definition of 
 a "space",  then spews out a "data structure" which is used for "past processing", that is for mathematical reasoning in that space. (In fact, when we say "let $M$ be a manifold", for example, we don't "have" that manifold,  only some properties of it, together with some very general abstract nonsense concerning 
 "legal" manipulations in the universe of "manifolds". All these can be liken with the image that we get 
 past the "front end" , in the sense that, like a real perceived image, we see it all, but we are incapable 
 of really enumerating and naming all that we see.)
 
 Even more, we may think that  the physical space can be understood, at some very fundamental level,   
 as the input of a "universal front end", and physical observers are "universal front ends". That is, maybe biology uses at a different scale an embodiment of a fundamental mechanism of the nature.

Thus, the  biologically inspired viewpoint is  that observers are like universal front ends looking at 
the same (but otherwise unknown) space. Interestingly this  may give a link between the problem of 
"local sign" in neuroscience and  the problem of understanding the properties of the physical space as emerging from some non-geometrical, more  fundamental structure, like a net, a foam, a graph... 

Indeed, in this physics research, one wants to obtain geometrical structure of the space (for example 
that locally, at the macroscopic scale,  it looks like $\mathbb{R}^n$) from a  non spatial like  structure, 
"seen from afar" (not unlike Gromov does with metric spaces). But in fact the brain does this all the time: 
from a class of intensive quantities (like the electric impulses sent by the neurons in the retina) the front 
end visual system reconstructs the space, literally as we see it. How it does it without "geometrical 
expertise" is called in neuroscience the problem of the "local sign" or of the "homunculus".  

\section{Simulating spaces}
 \label{sec1.2}
 
Let us consider then   some examples of spaces, like: the real world, a virtual world of a game, mathematical spaces as manifolds, fractals, symmetric spaces, groups, linear spaces ...
     
All these  spaces  may be characterized by the class of algebraic/differential computations which are possible, like: zoom into details, look from afar, describe velocities and perform other differential calculations needed for describing the physics of such a space, perform reflexions (as in symmetric spaces), linear combinations (as in linear spaces), do affine or projective geometry constructions and so on.
 
 Suppose that on a set $X$ (called "a space") there is an operation
$$ (x,y) \mapsto x \circ_{\varepsilon} y $$
which is dilation-like. Here $\varepsilon$ is a parameter belonging to a 
commutative group $\Gamma$, for simplicity let us take $\Gamma = (0,+\infty)$ 
with multiplication as the group operation. 

By "dilation-like" I mean the following: 
for any $x \in X$  the function 
$$ y \mapsto \delta^{x}_{\varepsilon} y \, = \, x \circ_{\varepsilon} y$$
behaves like a dilation in a vector space, that is 
\begin{enumerate}
\item[(a)] for any $\varepsilon, \mu \in (0,+\infty)$ we have 
$\displaystyle  \delta^{x}_{\varepsilon} \,  \delta^{x}_{\mu} \, = \,  \delta^{x}_{\varepsilon \mu}$ and $\displaystyle  \delta^{x}_{1} \, = \, id$; 
\item[(b)] the limit as $\varepsilon$ goes to $0$ of $\displaystyle 
 \delta^{x}_{\varepsilon} y$ is $x$, uniformly with respect to $x,y$. 
\end{enumerate} 

Then the dilation operation is the basic building block of both the 
algebraic structure of the space (operations in the tangent spaces) and the 
differential calculus in the space, as it will be explained further. This leads to he introduction 
of emergent algebras \cite{buligairq}. 

Something amazing happens if we take compositions of dilations, like this 
ones 
$$\Delta^{x}_{\varepsilon} (u,v) \, = \,  \delta_{\varepsilon^{-1}}^{\delta^{x}_{\varepsilon} u} \, \delta^{x}_{\varepsilon} u$$
$$\Sigma^{x}_{\varepsilon} (u,v) \, = \, \delta^{x}_{\varepsilon^{-1}} \, 
\delta_{\varepsilon}^{\delta^{x}_{\varepsilon} u} v $$
called the approximate difference, respectively approximate sum operations based 
at $x$. If we suppose that $\displaystyle(x,u,v) \mapsto \Delta^{x}_{\varepsilon} (u,v)$ and $\displaystyle(x,u,v) \mapsto \Sigma^{x}_{\varepsilon} (u,v)$ converge 
uniformly as $\varepsilon \rightarrow 0$ to $\displaystyle(x,u,v) \mapsto \Delta^{x} (u,v)$ and $\displaystyle(x,u,v) \mapsto \Sigma^{x} (u,v)$ then, out of apparently 
nothing, we get that $\displaystyle \Sigma^{x}$ is a group operation with $x$ as neutral element, which can be interpreted as the operation of vector addition 
in the "tangent space at $x$" (even if there is no properly defined such space). 

To convince you about this just look at the following example: $X = \mathbb{R~}^{n}$ and 
$$x \circ_{\varepsilon} y \,  = \, \delta^{x}_{\varepsilon} y \,  = \, x + \varepsilon (-x+y)$$ 
Then the approximate difference and sum operations based at $x$ have the expressions: 
\begin{align*}
\Delta_{\varepsilon}^{x}(u,v) 
&= x+\varepsilon(-x+u) + (-u+v)  & \rightarrow & & x -u+v \\
\Sigma_{\varepsilon}^{x}(u,v) 
&=  u+ \varepsilon(-u+x) + (-x+v)  & \rightarrow & & u -x+v 
\end{align*}

Moreover, with the same dilation operation we may define something resembling very 
much with differentiation. Take a function $f: X \rightarrow X$, then define 
$$ D_{\varepsilon} f (x) u \, = \, \delta_{\varepsilon^{-1}}^{f(x)} f \delta^{x}_{\varepsilon} (u)$$  
In the particular example used previously we get 
$$ D_{\varepsilon} f (x) u \, = \, 
f(x) + \frac{1}{\varepsilon} ( -f(x)+f(x+\varepsilon(-x+u)))$$
which shows that the limit as $\varepsilon$ goes to $0$ of $\displaystyle 
u \mapsto  D_{\varepsilon} f (x) u$ is a kind of differential.

Such computations are finite or virtually infinite "recipes", which can be implemented by some class of circuits made by very simple gates based on dilation 
operations (as in boolean computing, where transistors are universal gates for computing boolean functions).
     
(Computation in) a space is then described by emergent algebras \cite{buligairq},  which are inspired by the considerations about a formal calculus with binary decorated planar trees in relation with dilatation structures \cite{buligadil1}:
     
A - a class of transistor-like gates, with in/out ports labelled by points of the space and a internal state variable which can be interpreted as "scale". I propose dilations as such gates (basically these are idempotent right quasigroup operations).
     
B - a class of elementary circuits made of such gates (these are the "generators" of the space). The elementary circuits have the property that the output converges as the scale goes to zero, uniformly with respect to the input. 

C - a class of equivalence rules saying that some simple assemblies of elementary circuits have equivalent function, or saying that relations between those simple assemblies converge to relations of the space as the scale goes to zero.

Seen like this, "simulating a space" means: give a set of transistors (and maybe some non-emergent operations and relations), elementary circuits and relations which are sufficient to generate any interesting computation in this space.
     
For the moment I know how to simulate affine spaces \cite{buligadil2},  sub-riemannian  or Carnot-Carath\'eodory spaces \cite{buligadil3}, riemannian \cite{buligairq} or sub-riemannian symmetric spaces \cite{buligabraided}.

\section{Emergent algebras and decorated  braids diagrams}

Emergent algebras are idempotent right quasigroups (irqs) with a uniformity property. 

The  idempotent right quasigroups are related to algebraic structures appearing in knot theory.   J.C. Conway and G.C. Wraith, in their unpublished correspondence  from 1959,  used the name "wrack" for  a  self-distributive right quasigroup generated by a link diagram. Later, Fenn and Rourke \cite{fennrourke} proposed the name "rack" instead. Quandles are particular case of racks, namely self-distributive idempotent right quasigroups. They were introduced by  Joyce \cite{joyce}, as a distillation  of the  Reidemeister moves. 

The axioms of a (rack ;  quandle ; irq) correspond respectively to the (2,3 ; 1,2,3 ; 1,2) Reidemeister moves. That is why we shall use decorated braids diagrams in order to explain what energent algebras are.

The basic idea of racks and quandles is that these are algebraic operations related to the coloring 
of braids diagrams or, equivalently, of links diagrams.

In terms of braid colorings, with $X$ a set of colors, there are two binary operations on $S$ related 
to the coloring, as shown in the next figure. 

\centerline{\includegraphics[width=100mm]{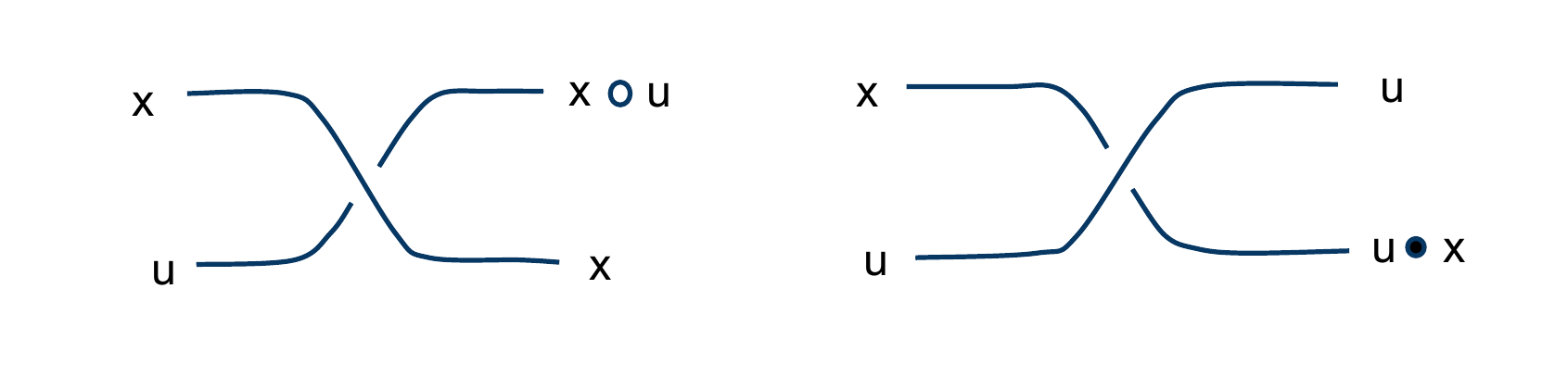}}

We have therefore a set $X$ endowed with two operations $\circ$ and $\bullet$ which satisfy a number of properties, such that the coloring of braids or links are compatible with the Reidemeister moves. 

The first Reidemeister move is depicted in the next figure, for braids and for links. 

\centerline{\includegraphics[width=120mm]{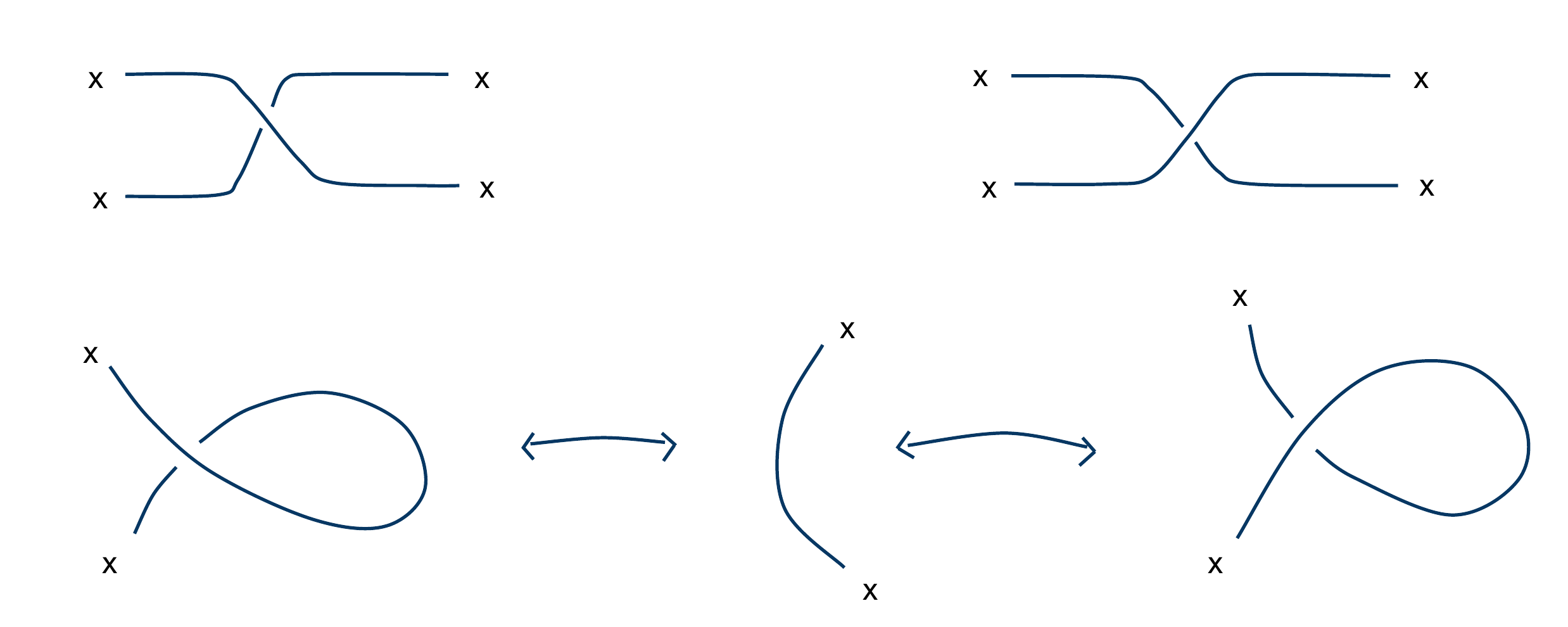}}

In algebraic terms we want the operations $\circ$ and $\bullet$ to be idempotent: 
$$x \circ x \, = \, x \bullet x \, = \, x$$
for all $x \in X$.

The second Reidemeister move, shown in the following picture, 

\centerline{\includegraphics[width=120mm]{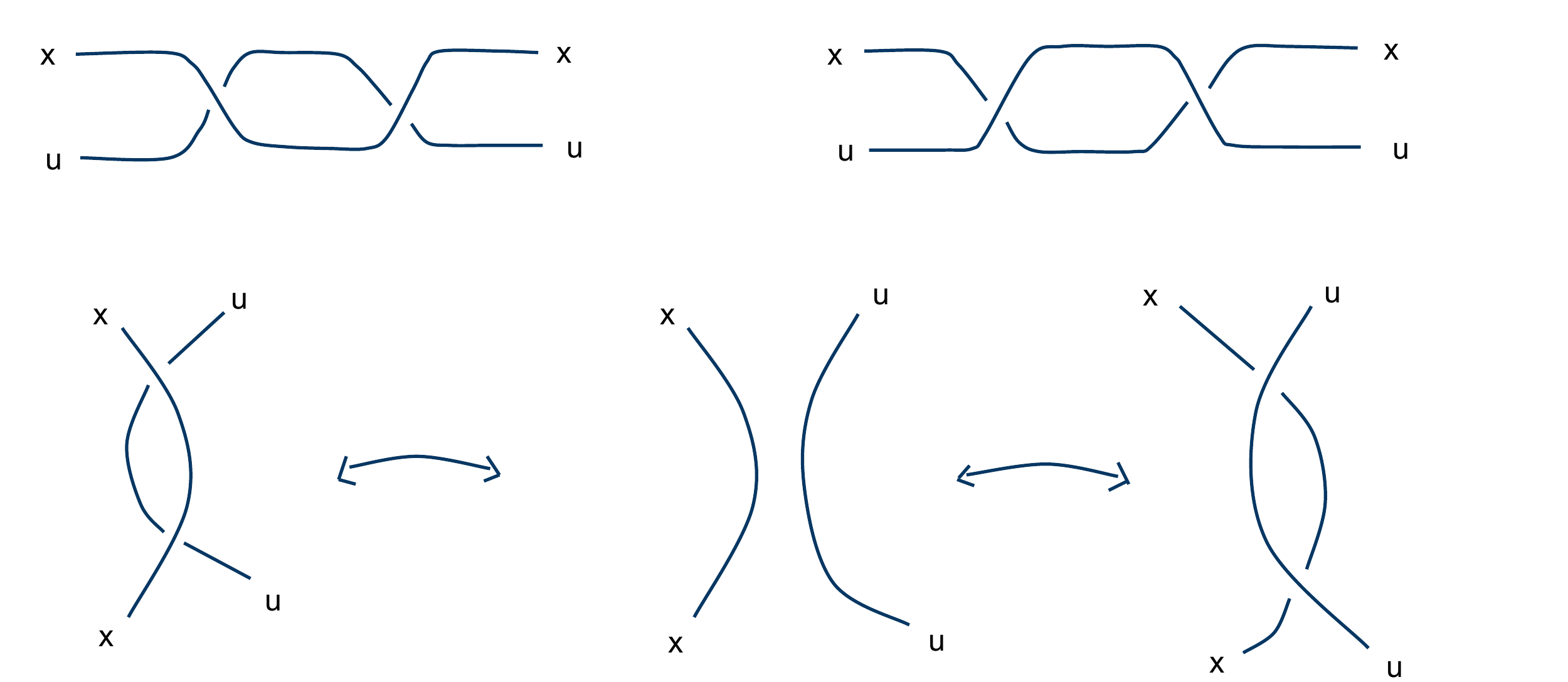}}

implies that we want $(X,\circ,\bullet)$ to be a right quasigroup, namely we want 
$$x \circ \left( x \bullet y\right) \, = \, x \bullet \left( x \circ y \right) \, = \, y$$
for all $x, y \in X$. This is the same as asking that for any $a$ and $b$ in $X$, the equation $a \circ x = b$ has a solution, which is unique, then denote the solution by $x = a \bullet b$.  
All in all, a set $(X,\circ,\bullet)$ which has the properties related to the first two Reidemeister 
moves is called an idempotent right quasigroup, or irq for short.

The third Reidemeister move: 

\centerline{\includegraphics[width=120mm]{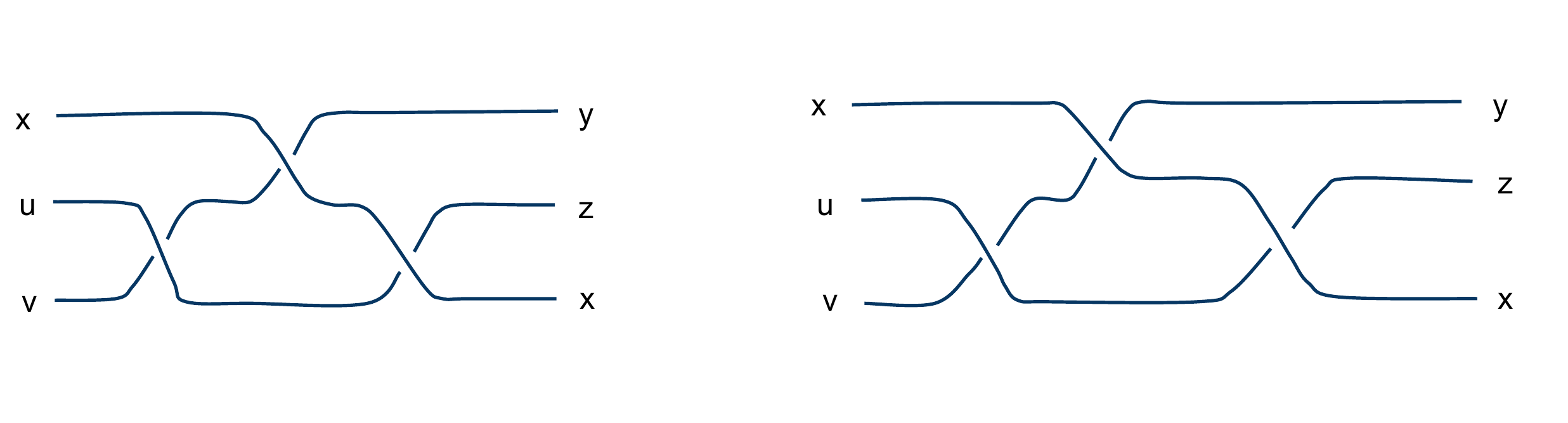}}

tells us, from the viewpoint  of the operations $\circ$ and $\bullet$, that we want $\circ$ (and therefore 
$\bullet$) to be self-distributive: 
$$x \circ \left( u \circ v \right) \, = \, \left( x \circ u \right) \circ \left( x \circ v \right)$$
Further we are not going  to ask for the operations of an emergent algebra (which is a one parameter 
family of irqs) to be self-distributive. Instead of this we shall want a weaker property to be satisfied. 

With this self-distributivity property, $(X,\circ,\bullet)$ is called a quandle. A well known quandle (therefore also an irq) 
is the Alexander quandle: consider $X = \mathbb{Z}[ \varepsilon, \varepsilon^{-1}]$ with the operations 
$$ x \circ y \, = \, x + \varepsilon \left( -x + y\right) \quad , \quad x \bullet y \, = \, x + \varepsilon^{-1} \left( -x + y\right)$$
So the operations in the Alexander quandle are dilations in euclidean spaces. 

Before passing to emergent algebras, let us write the  definition  of a irq.

\begin{definition} A right quasigroup is a set $X$ with a binary operation 
$\circ$ such that for each $a, b \in X$ there exists a unique $x \in X$ such that 
$a \, \circ\circ \, x \, = \, b$. We write the solution of this equation 
$x \, = \, a \, \bullet \, b$. 

 An idempotent right quasigroup (irq) is a  right quasigroup $(X,\circ)$ such that 
 for any $x \in X$ $x \, \circ \, x \, = \, x$. Equivalently, it can be seen as a 
  set $X$ endowed with two  operations $\circ$ and $\bullet$, which satisfy the following axioms: for any $x , y \in X$  
\begin{enumerate}
\item[(R1)] \hspace{2.cm} $\displaystyle x \, \circ \, x \, = \, x \, \bullet \, x \,  = \,  x$
\item[(R2)] \hspace{2.cm} $\displaystyle x \, \circ \, \left( x\, \bullet \,  y \right) \, = \, x \, \bullet \, \left( x\, \circ \,  y \right) \, = \, y$
\end{enumerate}
\label{defquasigroup}
\end{definition}

\subsection{Links with decorated crossings}

We could decorate not only the connected components of link diagrams or braids, but also the crossings. We shall use a commutative group $\Gamma$ for decorating the crossings, as shown in the next figure. 

\centerline{\includegraphics[width=120mm]{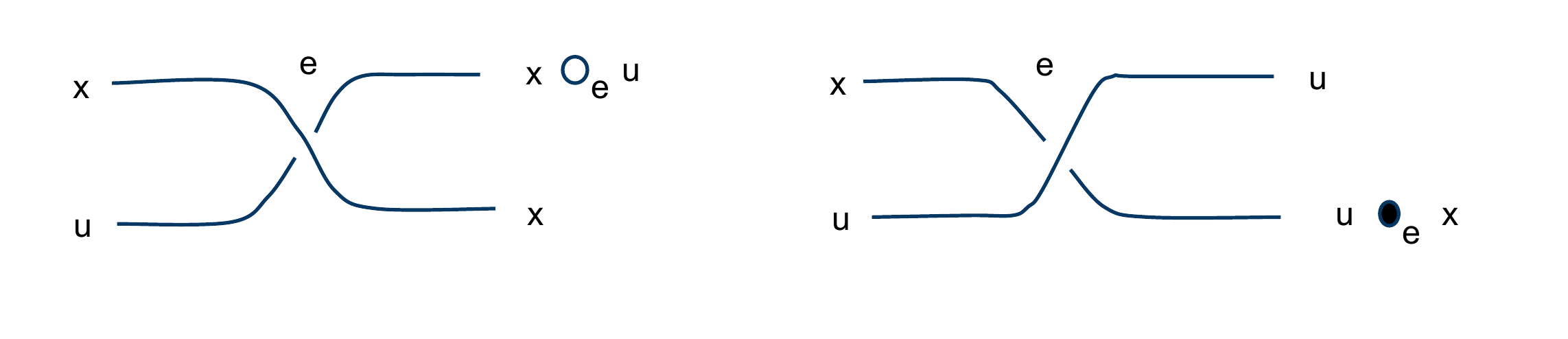}}

\paragraph{Notation convention.}Ê
In this figure, $e$ is an element of $\Gamma$. Further, in various formulae, we shall use greek letters 
for the elements of $\Gamma$,  like $\varepsilon$ instead of $e$, but in the following figures, as well, the elements of $\Gamma$ will be denoted with small latin letters.  This is just motivated by the graphic program available to me for producing the figures (google docs). 

The group $\Gamma$ which will be used mostly is $(0,+\infty)$ with multiplications. We shall explain things as if $\Gamma$ is this group. But keep in ming that other groups may be interesting, like 
the group of  nonzero complex numbers with multiplication, the group of integers with addition, or 
the direct product of any example from above with a finite commutative group.

We want the decoration to be such that for all $\varepsilon \in \Gamma$  the triples 
$\displaystyle (X, \circ_{\varepsilon}, \bullet_{\varepsilon})$ to be irqs, moreover we want that for 
any $x \in X$ the mapping 
$$ \varepsilon \in \Gamma \,  \mapsto \, x \circ_{\varepsilon} \, ( \cdot )$$
to be an action of $\Gamma$ on $X$. This is explained in the next two figures.

\centerline{\includegraphics[width=120mm]{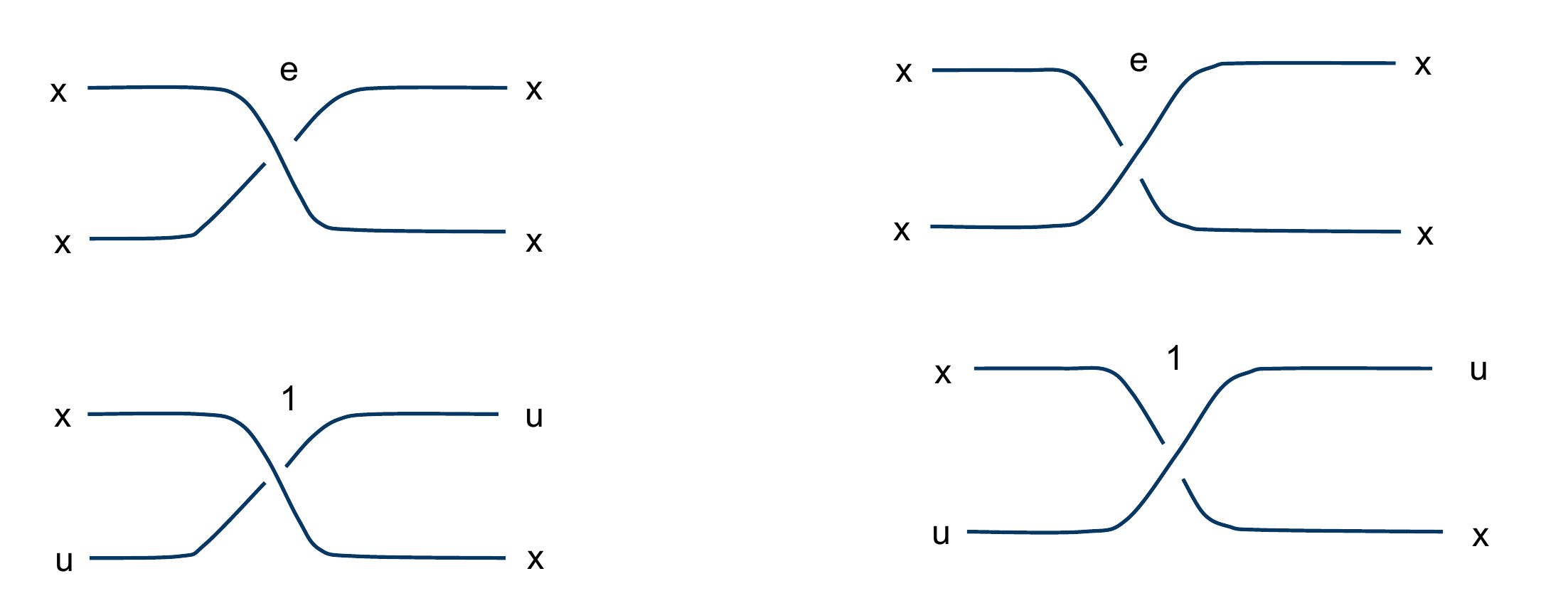}}

This figure tells us that the operations $\displaystyle \circ_{\varepsilon} ,  \bullet_{\varepsilon}$  are idempotent (first line of the figure) and moreover that $\displaystyle x  \circ_{1} u \, = u$ , 
$\displaystyle x  \bullet_{1} u \, = u$ for any $x, u \in X$ (second line of the figure). 
Here $1$ is the neutral element of $\Gamma$.

The next figure tells us that for equal inputs (from the left), the following decorated diagrams will have 
the same outputs (at the right). 

\centerline{\includegraphics[width=80mm]{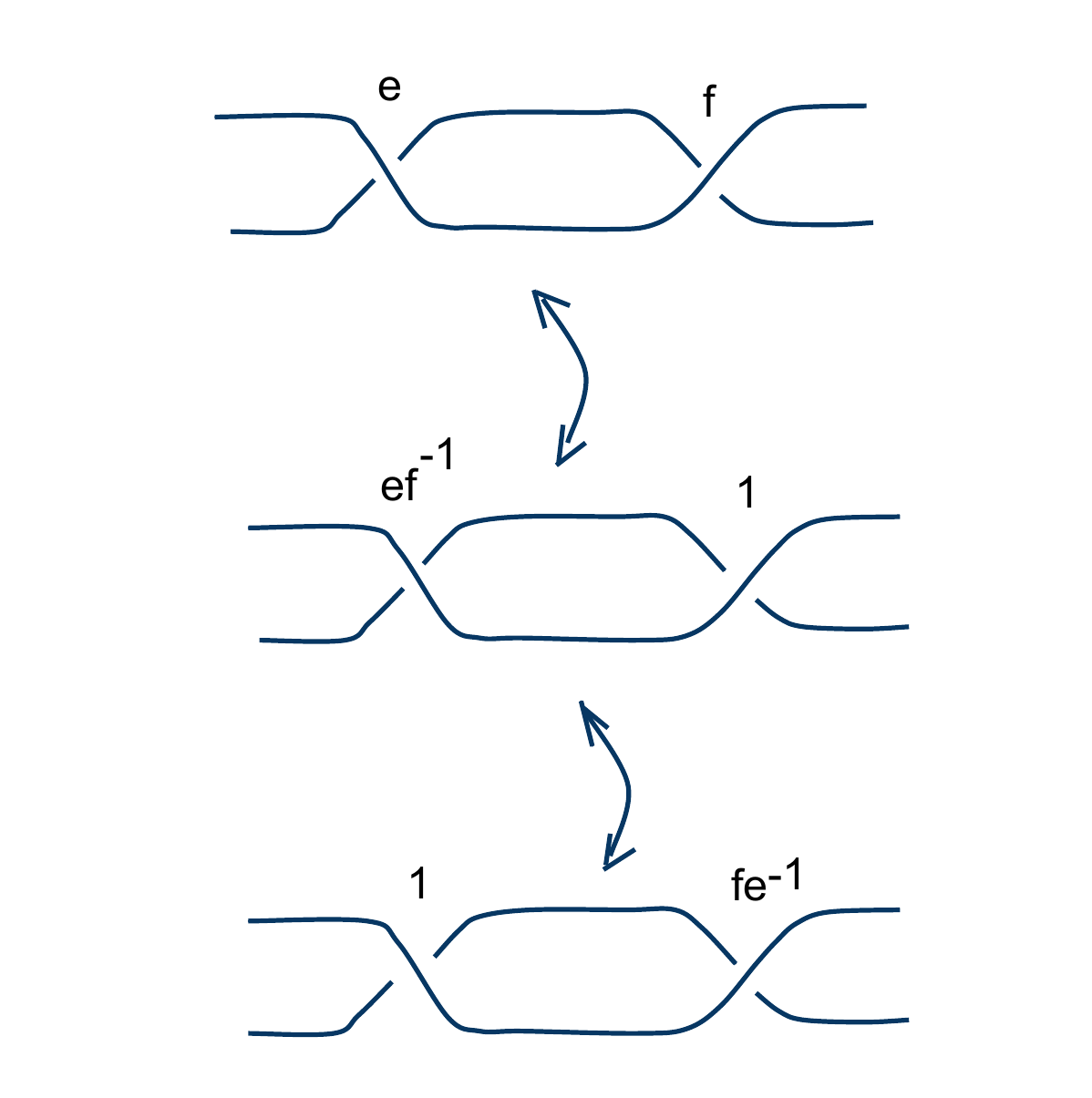}}

Taking into account the  second line of the previous figure,  parts of  this figure can be written as 
$$x  \circ_{\varepsilon} \left( x \bullet_{\mu} u \right) \, = \, x \circ_{\varepsilon \mu^{-1}} u$$
for any $x, u \in X$ and $\varepsilon, \mu \in \Gamma$. This implies  that the rule (R2), 
namely the second Reidemeister move, is true for $\displaystyle (X, \circ_{\varepsilon}, \bullet_{\varepsilon})$, and also that $\Gamma$ acts on $X$ in the sense explained before.

In \cite{buligairq} we introduced idempotent right quasigroups and then iterates
of the operations indexed by a parameter $\displaystyle k \in \mathbb{N}$. 
Here we used a general commutative group $\Gamma$, as in \cite{buligabraided}. 
The formal definition of a  $\Gamma$-irq is given further.

\begin{definition}
Let $\Gamma$ be a commutative group. A $\Gamma$-idempotent right quasigroup 
is a set $X$ with a function $\displaystyle \varepsilon \in \Gamma \mapsto 
\circ_{\varepsilon}$ such that for any $\varepsilon \in \Gamma$ the pair $\displaystyle (X, \circ_{\varepsilon})$ is a irq
and moreover for any $\varepsilon, \mu \in \Gamma$ and any $x, y \in X$ we have 
$$x \, \circ_{\varepsilon} \, \left( x \, \circ_{\mu} \, y \right) \, = \, 
x \, \circ_{\varepsilon \mu} \, y$$
\label{defgammairq}
\end{definition}

The approximate difference gate, introduced in the next definition, is central to the 
reasonings of this paper .

\begin{definition}
For any $\varepsilon \in \Gamma$, the $\varepsilon$-approximate difference gate, described by the next figure as the gate $\displaystyle DIF_{\varepsilon}(x,u,v) \, = \, (x, x\circ_{\varepsilon}u, \Delta_{\varepsilon}^{x}(u,v))$.
\label{defdifapprox} 
 Here   $\displaystyle  \Delta^{x}_{\varepsilon} (u,v)$ is a construct made from operations $\displaystyle \circ_{\varepsilon}$, $\displaystyle \bullet_{\varepsilon}$. From the figure we 
can compute $\displaystyle  \Delta^{x}_{\varepsilon} (u,v)$ as 
$$\displaystyle  \Delta^{x}_{\varepsilon} (u,v) \, = \, \left( x \circ_{\varepsilon} u \right) \, 
\bullet_{\varepsilon} \, \left( x \circ_{\varepsilon} v \right)$$
\end{definition}
\centerline{\includegraphics[width=120mm]{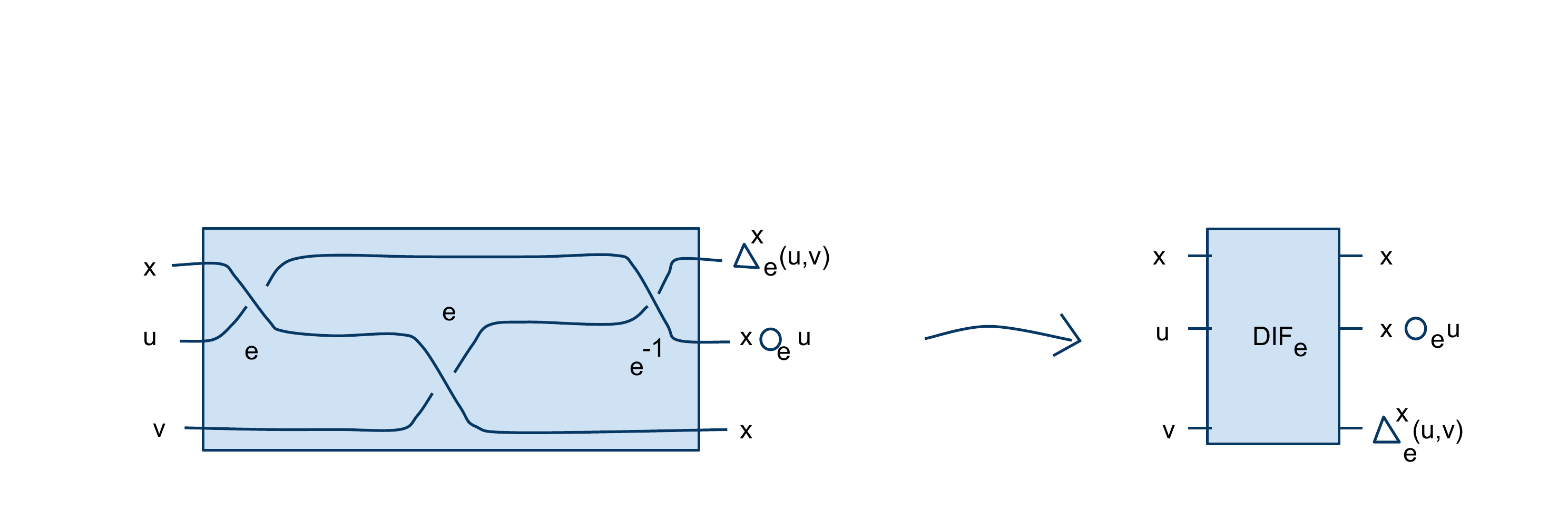}}

There are two more constructs made from the operations of a $\Gamma$-irq, namely the approximate 
sum gate and the approximate inverse gate. The approximate sum gate is figured below. 

\centerline{\includegraphics[width=120mm]{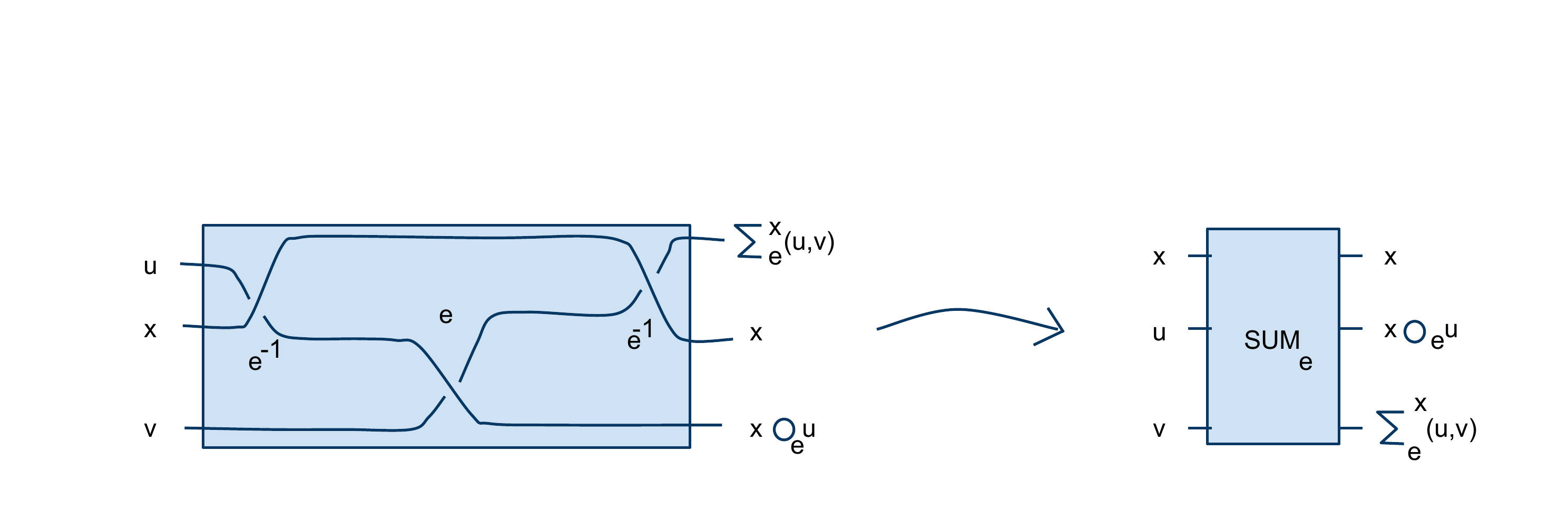}}

The inverse gate will be useful, too. It is shown in the next figure. 

\centerline{\includegraphics[width=80mm]{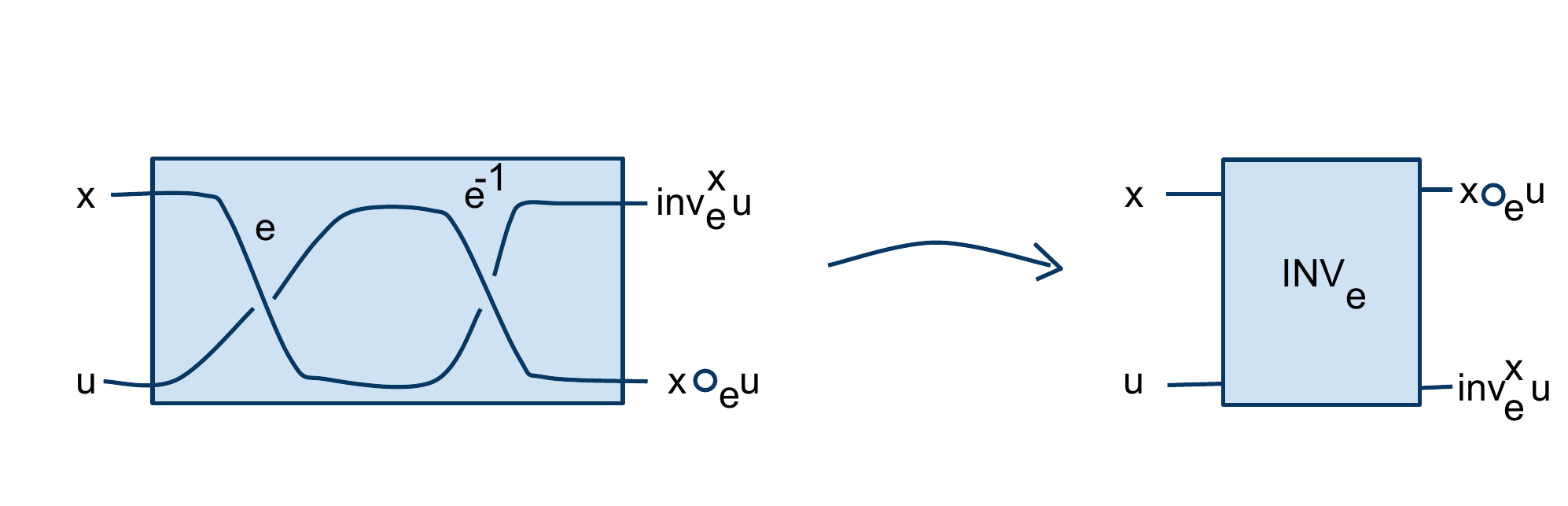}}

The geometric meaning of $\displaystyle \Delta^{x}_{\varepsilon}(u,v)$ is that it is indeed a kind of approximate difference between the vectors $\displaystyle \vec{xu}$ and $\vec{xv}$, by means of a generalization of the 
parallelogram law of vector addition.  This is  shown in the following figure, where 
straight lines have been replaced by slightly curved ones in order to suggest that this construction has meaning in settings far more general than euclidean spaces, like in Carnot-Caratheodory or 
sub-riemannian geometry, as shown in \cite{buligadil1}, or generalized (noncommutative) affine geometry 
 \cite{buligadil2}, for length metric spaces with dilations \cite{buligadil3} or even for normed groupoids 
 \cite{buligagr}. 
\centerline{\includegraphics[width=120mm]{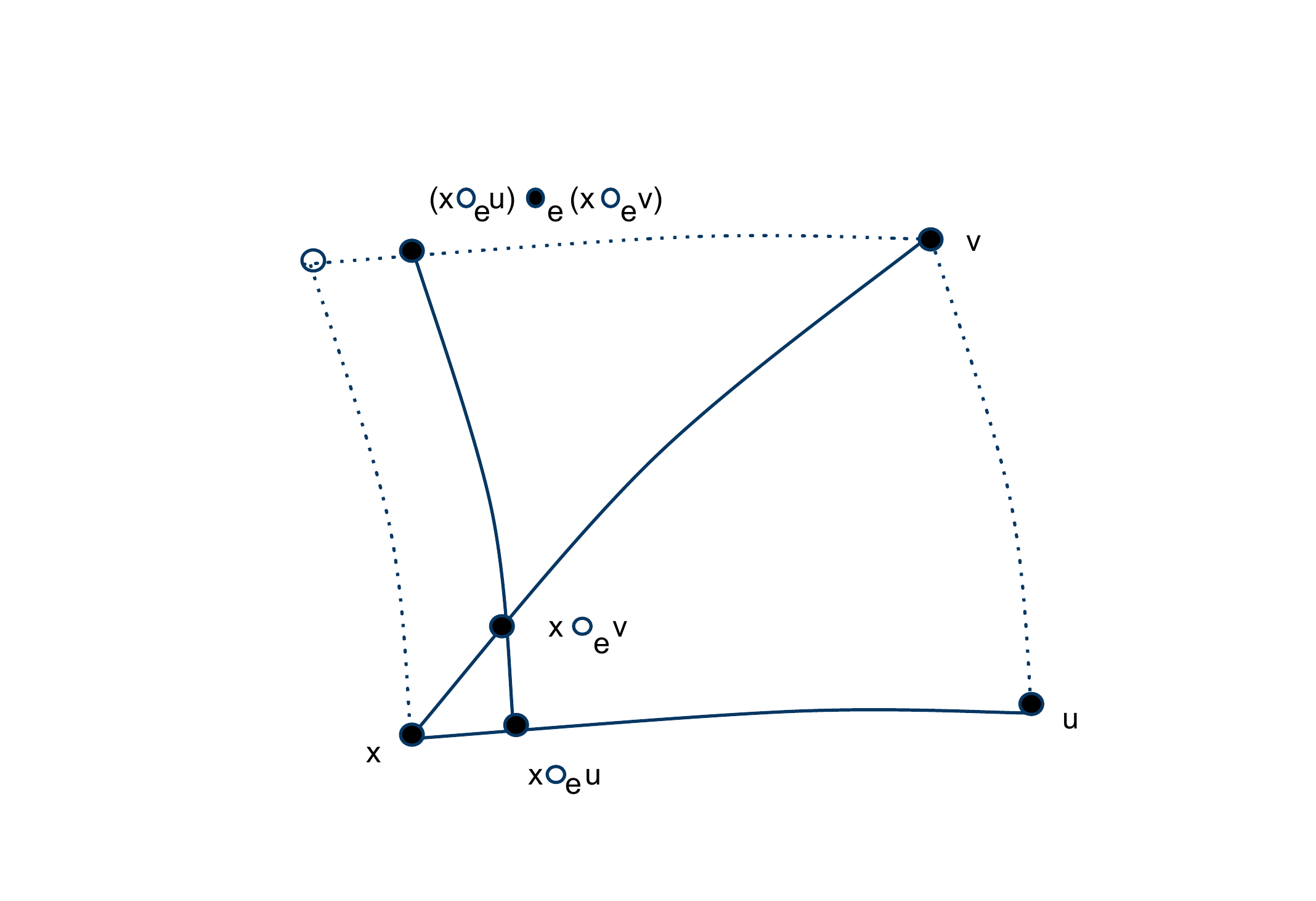}}

In order to better understand these constructs with dilations, in \cite{buligadil1} was proposed a decorated binary trees notation, which makes further computations easier.  

\subsection{Decorated planar binary trees} 
The idea is to see the outputs of these gates as particular trees, constructed from the following two 
generators: 

\centerline{\includegraphics[width=120mm]{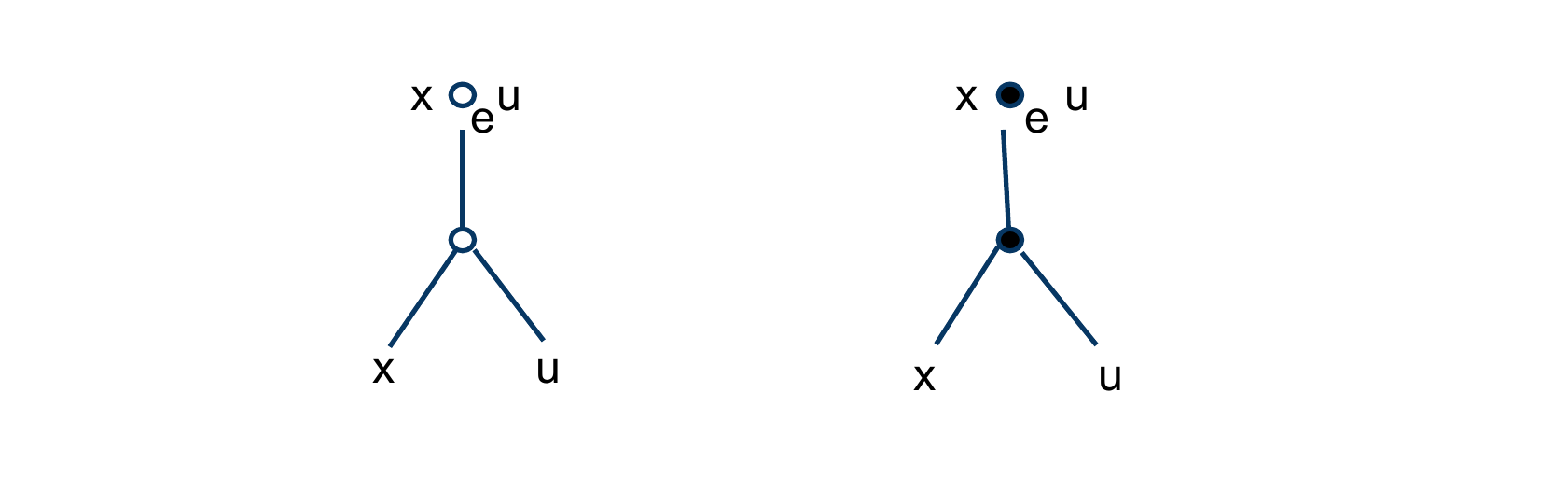}}

In the language of trees,  Reidemeister rules 1 and 2 from 
knot theory correspond to the rules (R1), (R2) described below. 

\centerline{\includegraphics[width=120mm]{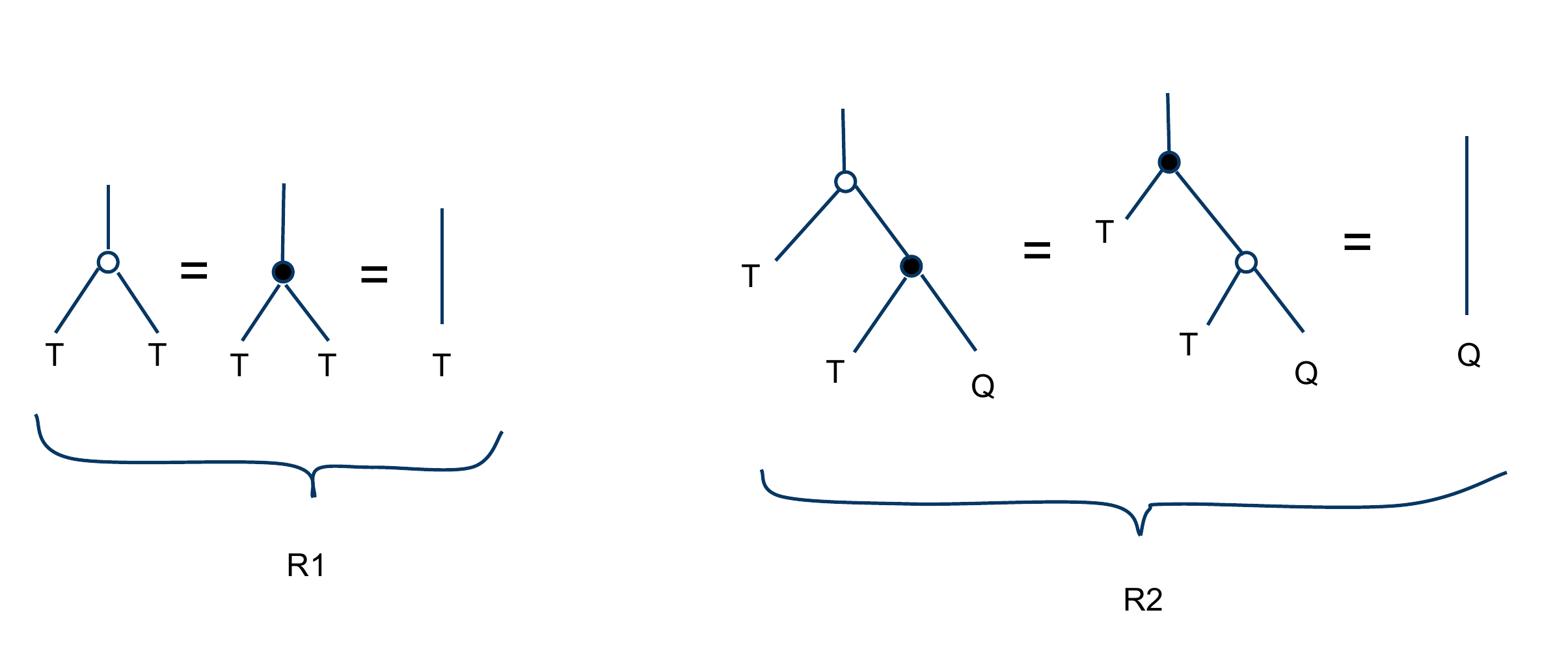}}

The relevant outputs of the previously introduced gates, namely the approximate difference, sum and 
inverse functions, are described in the next definition. 

\begin{definition}
We define the difference, sum and inverse trees  given by:

\centerline{\includegraphics[width=120mm]{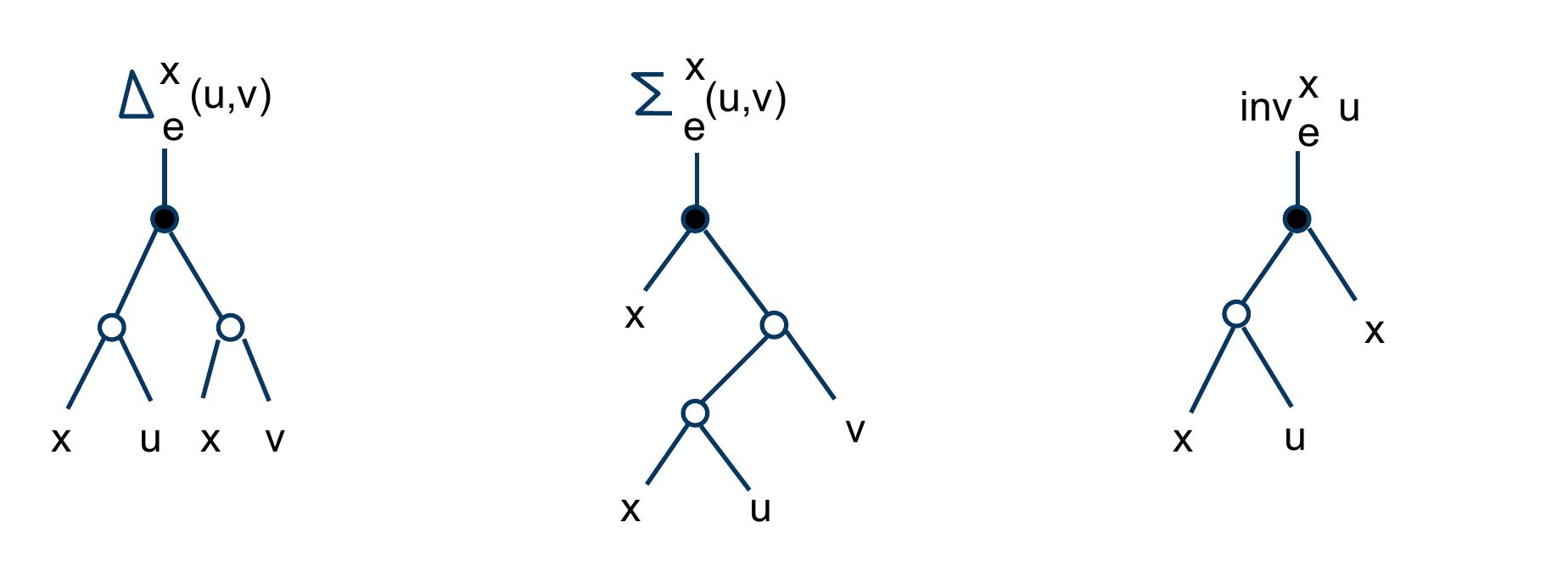}}

\label{defsumdif}
\end{definition}

\begin{proposition}
Let $\displaystyle (X,\circ_{\varepsilon})_{\varepsilon \in \Gamma}$ be
a $\Gamma$-irq. Then we have the relations: 
\begin{enumerate}
\item[(a)] (difference is the inverse of sum) $\displaystyle \Delta^{x}_{\varepsilon}(u, \Sigma^{x}_{\varepsilon}(u,v)) \, = \, v$ 
\item[(b)] (sum is the inverce of difference)$\displaystyle \Sigma^{x}_{\varepsilon}(u, \Delta^{x_{\varepsilon}}(u,v)) \, = \, v$ 
\item[(c)] (difference  approximately equals  the sum of the inverse)  $\displaystyle \Delta^{x}_{\varepsilon}(u, v) \, = \, \Sigma^{x \circ_{\varepsilon} u}_{\varepsilon} 
(inv_{\varepsilon}^{x} u , v)$
\item[(d)] (inverse operation is approximatively an involution)  $\displaystyle inv_{\varepsilon}^{x\circ u} \, inv_{\varepsilon}^{x} \, u  \, = \, u $
\item[(e)] (approximate associativity of the sum)  $\displaystyle \Sigma^{x}_{\varepsilon}(u, \Sigma^{x\circ_{\varepsilon}u}_{\varepsilon}
(v , w)) \, = \, \Sigma^{x}_{\varepsilon}(\Sigma^{x}_{\varepsilon}(u,v), w) $
\item[(f)] $\displaystyle  inv^{x}_{\varepsilon} \, u \, = \,  \Delta^{x}_{\varepsilon}( u , x)$
\item[(g)] (neutral element at right)  $\displaystyle  \Sigma^{x}_{\varepsilon} (x, u) \, = \,  u $
\end{enumerate}
\label{pplay}
\end{proposition}

We shall use the tree formalism to prove  some 
of these relations.   For complete proofs see \cite{buligadil1}. 

For example, in order to prove  (b) we do the following calculus:

\centerline{\includegraphics[width=120mm]{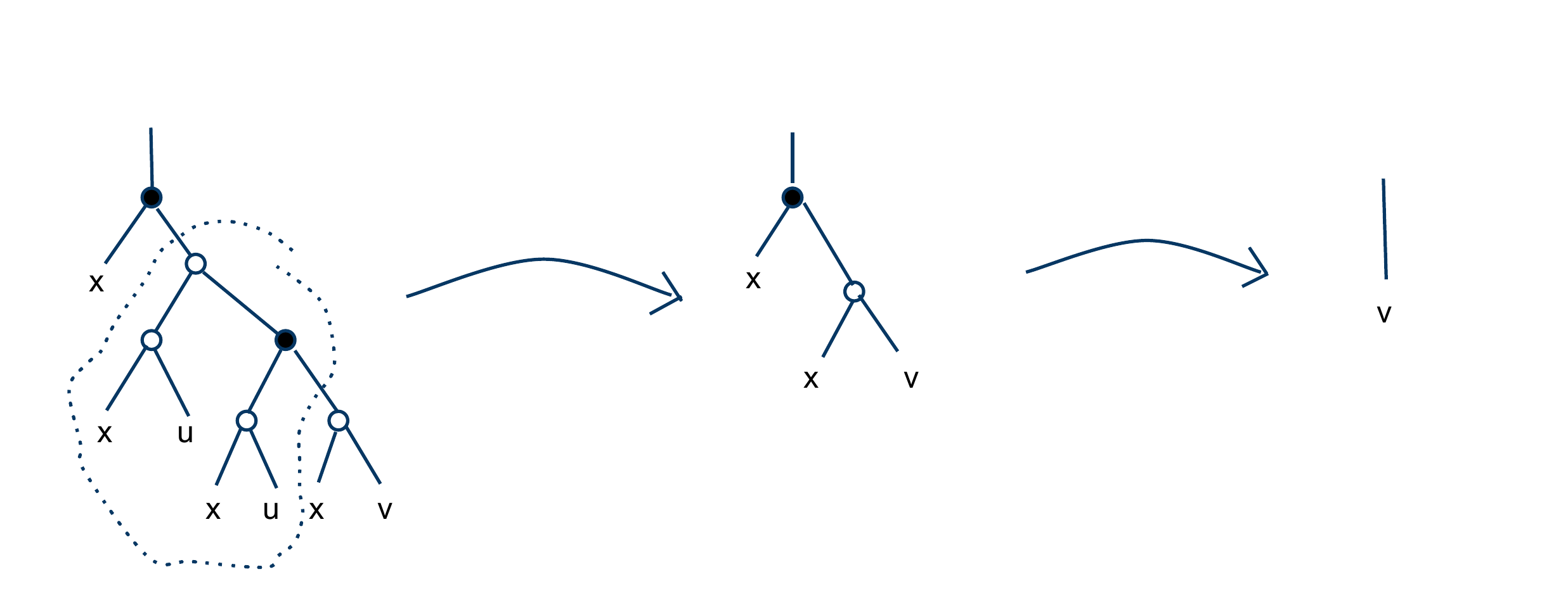}}

The relation (c) is obtained from: 

\centerline{\includegraphics[width=100mm]{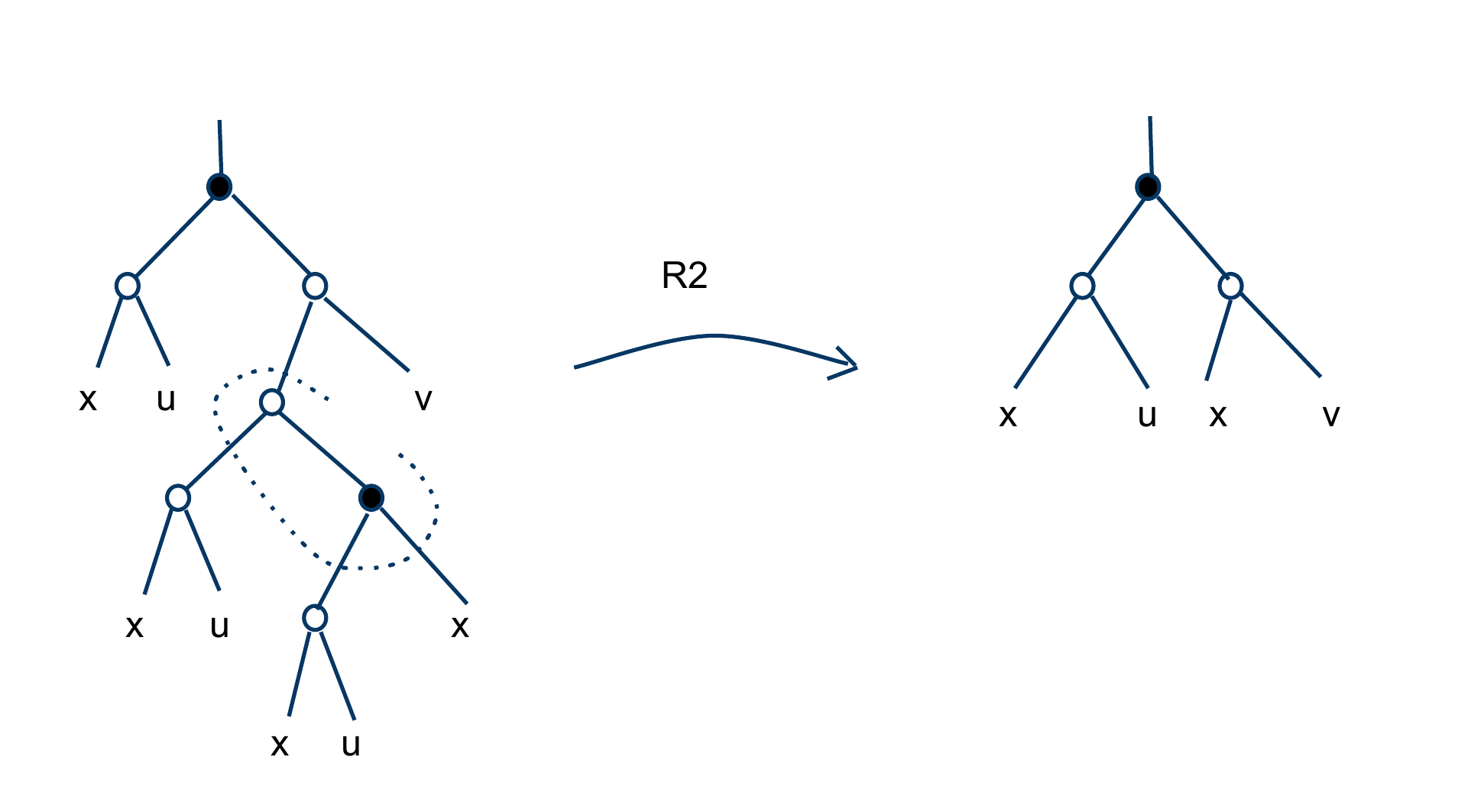}}

Relation (e) (which is  a kind of associativity relation) is obtained from: 

\centerline{\includegraphics[width=120mm]{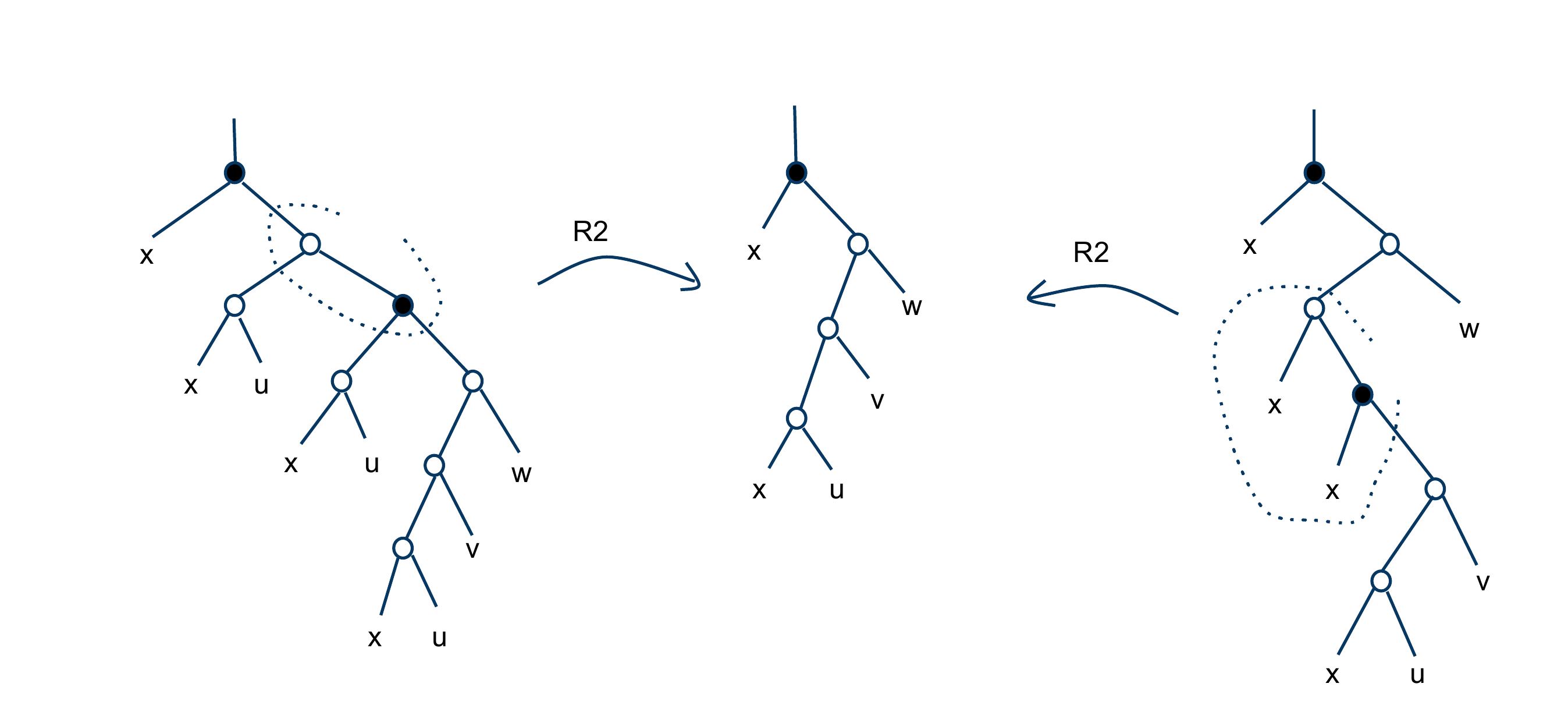}} 

Finally, for proving  relation (g) we use also the rule (R1). 

\centerline{\includegraphics[width=80mm]{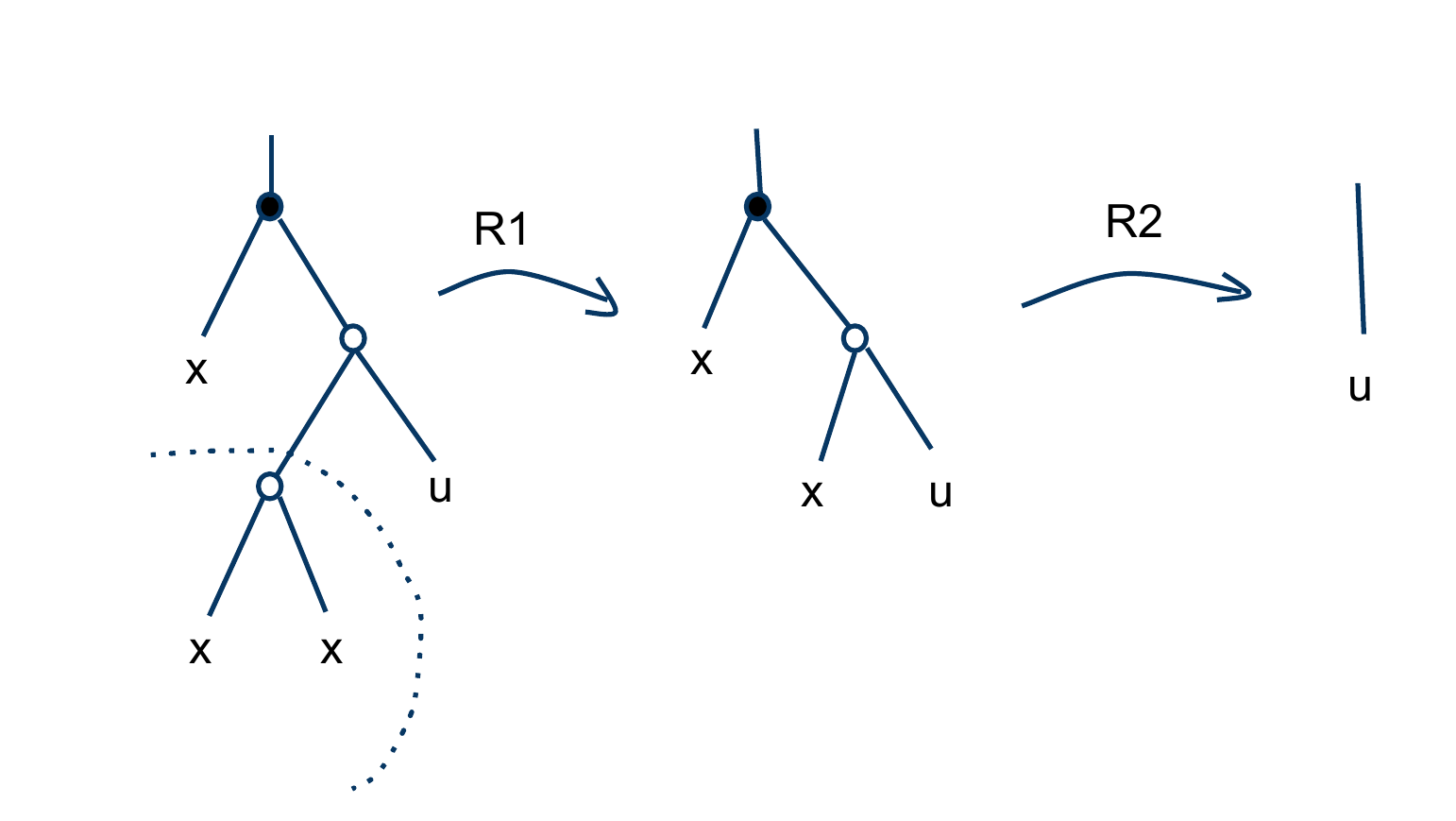}}

\subsection{Arbitrarily good approximations}

We arrive now at the most important axioms of emergent algebras. An emergent algebra is 
a $\Gamma$-irq with two properties that taken together are weaker than the Reidemeister 
move 3 law, see section \ref{r3weaker}. 

We shall need to put on $\Gamma$ a supplementary structure, namely an absolute. 
In order not to complicate things, let us take $\Gamma = (0,+\infty)$, thus we have no difficulty 
in saying what we mean by $\Gamma \ni \varepsilon  \rightarrow 0$. (More general, we might 
take $\mid \cdot \mid : \Gamma \rightarrow (0,+\infty)$ a group morphism and define 
"$\varepsilon \rightarrow 0$" by $\mid \varepsilon \mid \rightarrow 0$.)

We shall now suppose that $X$ is a topological space endowed with an uniform structure.  
Uniform structures are stronger than topologies, technically they are filters of the diagonal 
$\left\{ (x,x) \mbox{Ê: } x \in X \right\} \subset X \times X$ with certain poperties, which allow 
us to say that something is happening " uniformly with respect to $x, y, .... \in X$".  

\begin{remark}
The uniformity 
is essential for further reasonings, see theorem \ref{mainthm}.
\end{remark}

The two supplementary axioms that we need  are saying that something happens when  
$\varepsilon \rightarrow 0$, uniformly with respect to the entries. The first axiom is the following. 

\centerline{\includegraphics[width=120mm]{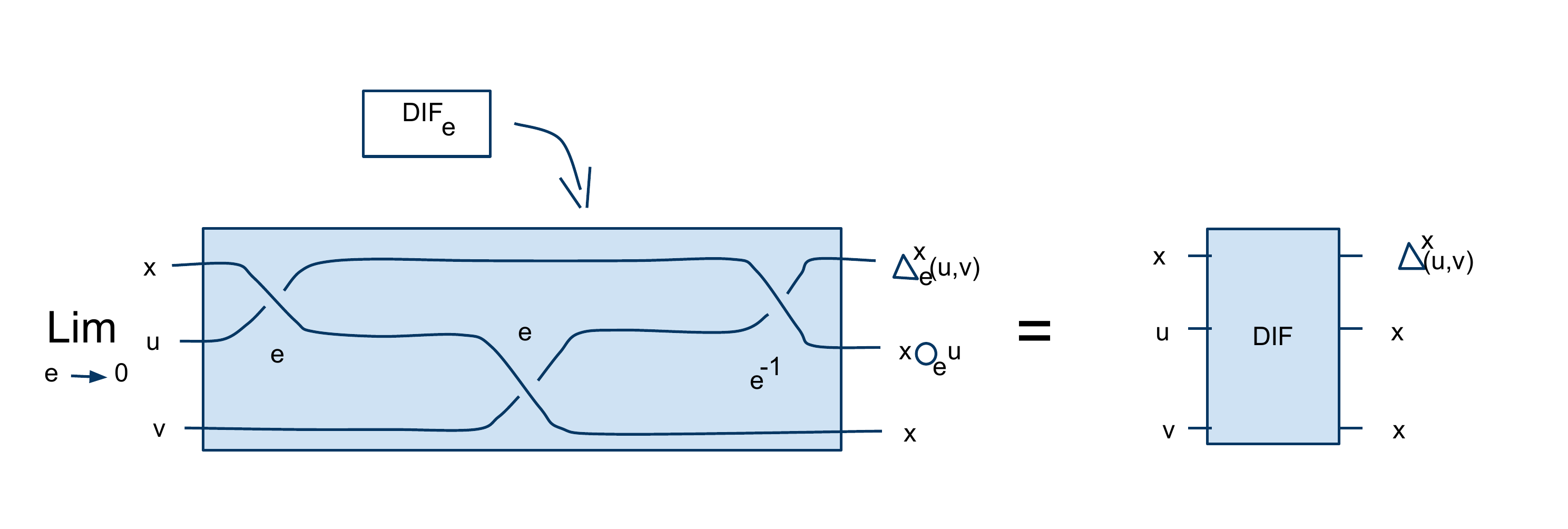}}

Therefore, as $\varepsilon \rightarrow 0$, we have 
$$\lim_{\varepsilon \rightarrow 0} x \circ_{\varepsilon} u \, = \, x$$
uniformly with respect to $x , u$ in a arbitrary compact subset of   $X$. Likewise, 
$$  \lim_{\varepsilon \rightarrow 0} \Delta^{x}_{\varepsilon} (u,v) \, = \, \Delta^{x}(u,v)$$ 
uniformly with respect to $x, u , v$ in a  arbitrary compact subset of $X$.

\begin{remark}
At this stage, you may want to read again section \ref{sec1.2}, in order to understand where 
this is going. 
\end{remark}

The axiom tells us that gates $\displaystyle DIF_{\varepsilon}$ approximate 
arbitrarily well (and uniformly with respect to the inputs) the gate 
$$DIF(x,u,v) \, = \, (x, x, \Delta^{x}(u,v))$$

\begin{remark}
The idea of gates which approximate arbitrarily well a given gate is present in Deutsch formulation 
of quantum computing \cite{deut}. 
\end{remark}

Finally, the second axiom is the following. 

\centerline{\includegraphics[width=120mm]{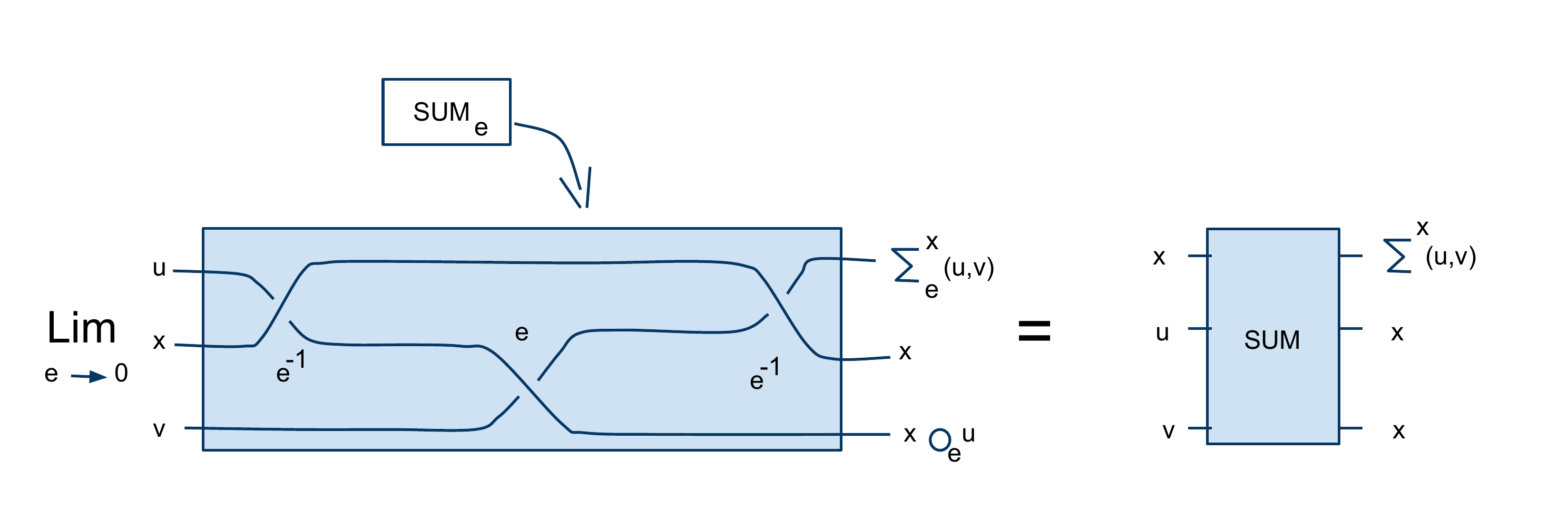}}
The new thing stated by this axiom is that 
$$  \lim_{\varepsilon \rightarrow 0} \Sigma^{x}_{\varepsilon} (u,v) \, = \, \Sigma^{x}(u,v)$$ 
uniformly with respect to $x, u , v$ in a  arbitrary compact subset of $X$. Here   $\displaystyle  \Sigma^{x}_{\varepsilon} (u,v)$ can be computed from the figure,  as 
$$\displaystyle  \Sigma^{x}_{\varepsilon} (u,v) \, = \, x \bullet_{\varepsilon} \left( \left( 
x \circ_{\varepsilon}  u \right) \, \circ_{\varepsilon} v   \right)$$

\begin{remark}
If $X$ is locally compact then the second axiom is not needed, because it can be deduced from the first. 
\end{remark}

We are now going to give precise definitions. See \cite{buligabraided} for all details concerning absolutes and uniformities, or take a look at the Appendix. 

\begin{definition}
A $\Gamma$-uniform irq, or emergent algebra  $(X, \circ, \bullet)$ is a separable uniform  
space $X$ which is also a $\Gamma$-irq, with continuous operations,  such that: 
\begin{enumerate}
\item[(C)] the operation $\circ$ is compactly contractive: for each compact set 
$K \subset X$ and open set $U \subset X$, with $x \in U$, there is an open set 
$\displaystyle A(K,U) \subset \Gamma$ with $\mu(A) = 1$ for any $\mu \in 
Abs(\Gamma)$ and for any $u \in K$ and
$\varepsilon \in A(K,U)$, we have $\displaystyle x \circ_{\varepsilon} u \in U$; 
\item[(D)] the following limits exist for any $\mu \in Abs(\Gamma)$ 
$$ \lim_{\varepsilon \rightarrow \mu} \Delta_{\varepsilon}^{x}(u,v)  \, = \, 
\Delta^{x}(u,v) \quad , \quad \lim_{\varepsilon \rightarrow \mu} \Sigma_{\varepsilon}^{x}(u,v) 
 \,  = \, \Sigma^{x}(u,v) $$
and are uniform with respect to $x, u, v$ in a compact set. 
\end{enumerate}
\label{deftop}
\end{definition}

The main property of a uniform irq is the following. It is a consequence of 
relations from proposition \ref{pplay}. 

\begin{theorem}
Let $(X, \circ, \bullet)$ be a uniform irq. Then for any $x \in X$  the operation 
$\displaystyle (u,v) \mapsto  \Sigma^{x}(u,v) $ gives $X$ the structure of
a conical group with the dilation $\displaystyle u \mapsto x \circ_{\varepsilon} u$.
\label{mainthm}
\end{theorem}

\paragraph{Sketch of the proof.} Pass to the limit in the relations from proposition \ref{pplay}.  We can do this exactly because of the uniformity assumptions. We therefore have a series of algebraic relations which can be used to get the conclusion.

\subsection{Reidemeister 3 move "emerges" from the first two moves}
\label{r3weaker}

In this section I explain why, from the point of view of  emergent algebras, the last two axioms 
imply  the Reidemeister 3 move. For this we have to understand better the content of the theorem 
\ref{mainthm}.

In the following figure we see a particular tangle, called relative dilation with respect to $(x,\varepsilon) \in X \times \Gamma$.  

\centerline{\includegraphics[width=120mm]{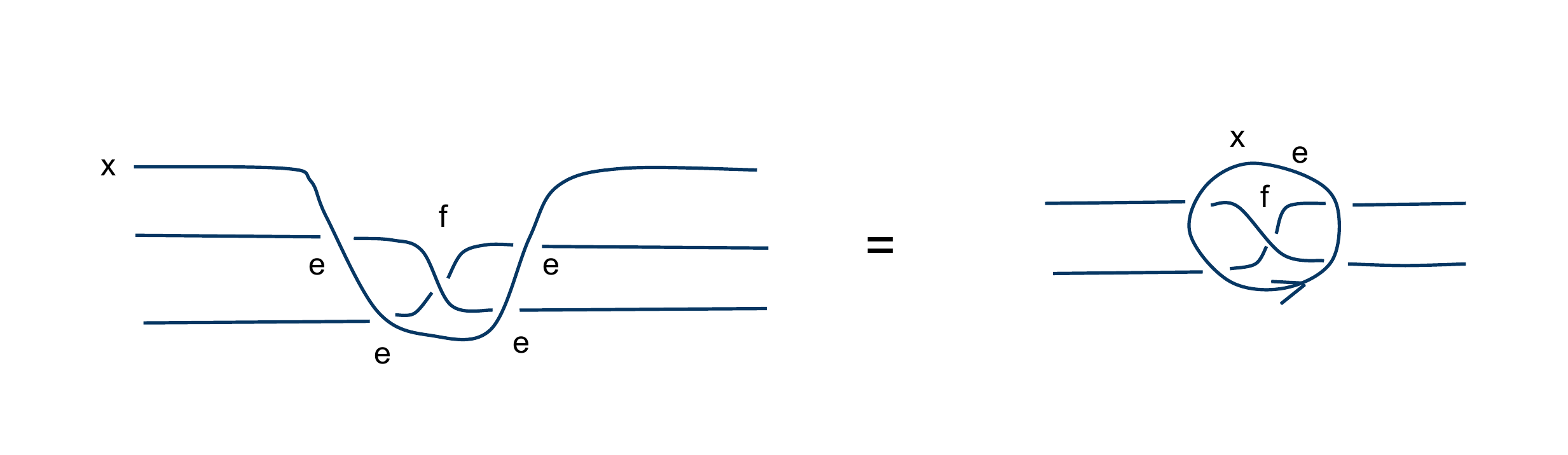}}

Seen as an operation, the relative dilation has the following expression, for fixed  $(x,\varepsilon) \in X \times \Gamma$:  for every $\lambda \in \Gamma$ 
$$ (u,v) \in X \times X  \, \mapsto \,   u \circ^{x,\varepsilon}_{\lambda} v \, = \, x \bullet_{\varepsilon} ((x\circ_{\varepsilon}u)\circ_{\lambda} (x \circ_{\varepsilon}v))$$
It is clear that $\displaystyle (X, \circ^{x,\varepsilon})$ is a $\Gamma$-irq.  Moreover, 
$(X,\circ)$ is self-distributive if and only if for any $(x,\varepsilon) \in X \times \Gamma$ we have 
$$(X, \circ^{x,\varepsilon}) \, = \, (X,\circ)$$
Self-distributivity is related to the Reidemeister move (R3), so the previous equality characterizes 
the move (R3) from the viewpoint of emergent algebras. 

An interesting consequence of the theorem \ref{mainthm}Ê   is the following. We shall suppose that 
$\Gamma = (0,\infty)$ and the absolute is $0$. 

\begin{corollary}
In the hypothesis of theorem \ref{mainthm},  $\displaystyle (X, \circ^{x,\varepsilon})$ converges as 
$\varepsilon \rightarrow 0$ to a self-distributive emergent algebra, in the sense
$$\lim_{\varepsilon \rightarrow 0} u \circ^{x,\varepsilon}_{\lambda} v \, = \, 
\Sigma^{x}(u, \cdot) (x \circ_{\lambda} \cdot) \Delta^{x}(u,v)$$
\end{corollary}

The fact the operation from the right hand side is self-distributive is a consequence of the characterization of conical groups as self-distributive emergent algebras. 

For the proof of the corollary let us decompose the relative dilation into three parts, as 
shown here. 

\centerline{\includegraphics[width=120mm]{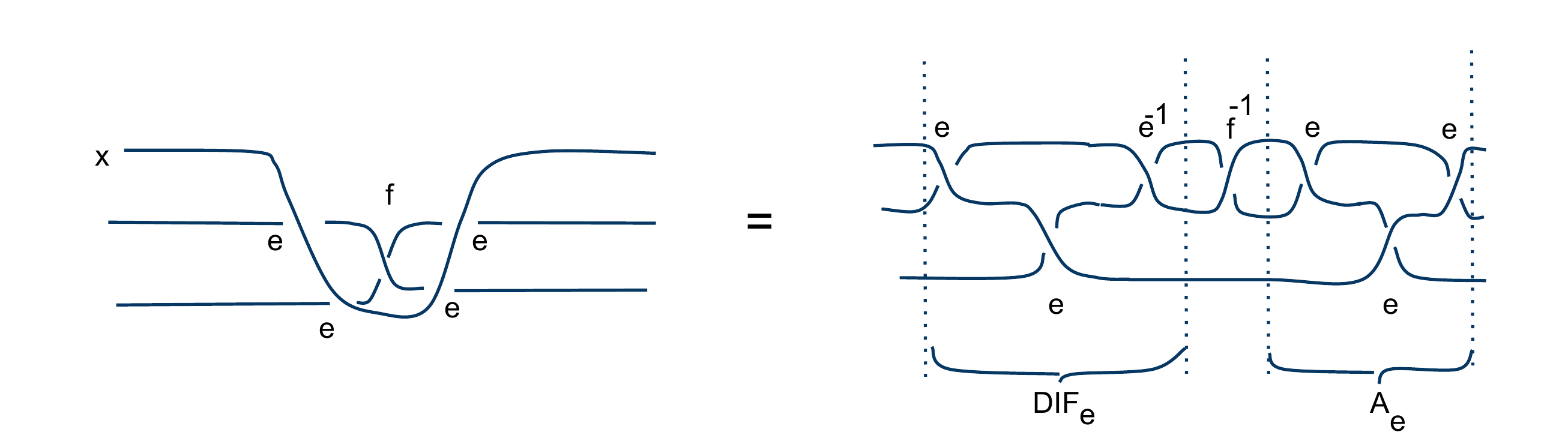}}

The first part is a difference gate. The second part is similar to a sum gate, in the sense explained 
in the next figure:

\centerline{\includegraphics[width=120mm]{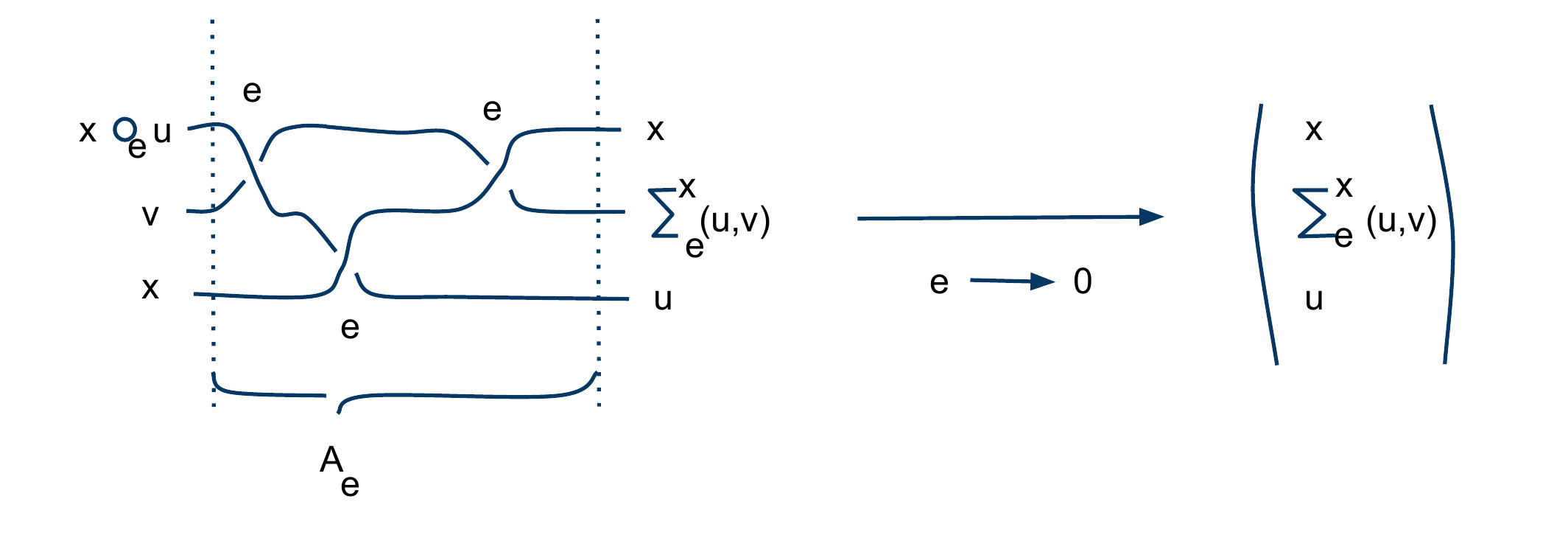}}

Therefore as $e \rightarrow 0$ we obtained the desired conclusion.

\section{A tentative interpretation of the front end by emergent algebras}

In Meredith and  Snyder \cite{meredith} the authors propose an encoding of knots in the $\pi$-calculus 
\cite{pi}, such that knots are ambient isotopic if and only if their encodings are weakly bisimilar. In simpler words they attach to any knot diagram a process such that the processes associated to the knot diagram 
before and after one of the Reidemeister moves are related by a bisimulation.

In this section I am going to use the idea that the elements of the space $X$, which I use 
to decorate links, can be seen as names in $\pi$-calculus, which have both roles of variables 
and communication channels. In the front end interpretation $x \in X$ is either a name 
for a point of the visual system (thus a name of a channel) or a variable (for example intensity of the 
light field) communicated through that channel. 

With this interpretation, the operation $\displaystyle x \circ_{\varepsilon} y$ has two meanings, namely 
if $x, y$ are names of channels then a dilation gate serves to define another channel, while if 
$x,y$ refer to the variables communicated through channels then the dilation gate computes a 
weighted average of the variables and communicates it to the channel just defined. 

Therefore the points of the visual system are encoded by a set $X$ of names, either of channels, or of 
variables, which has the structure of an emergent algebra. 

The interpretation of a cortical hypercolumn in this (sketch of a) formalism over $x \in X$ is the 
functor which associates to any decorated link diagram in braids notation the relative to $x$ diagram, as 
figured here. 

\centerline{\includegraphics[width=120mm]{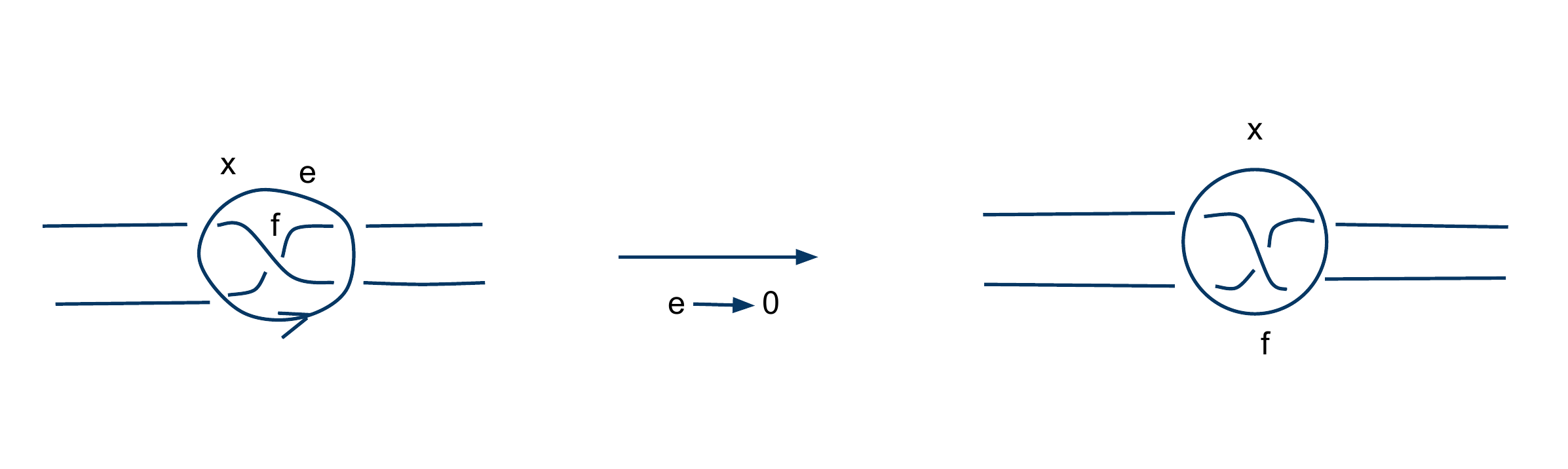}}

Here we have a braid diagram which is "encircled" with a curve decorated with $x$, such that 
all new crossings which appear are decorated with $\varepsilon$.

This is indeed a functor, in the sense that it is a morphism of decorated braids . The operation of encircling the whole link diagram is equivalent, by rules (R1) and (R2) only, to the operation of encircling 
of each crossing. 

\centerline{\includegraphics[width=70mm]{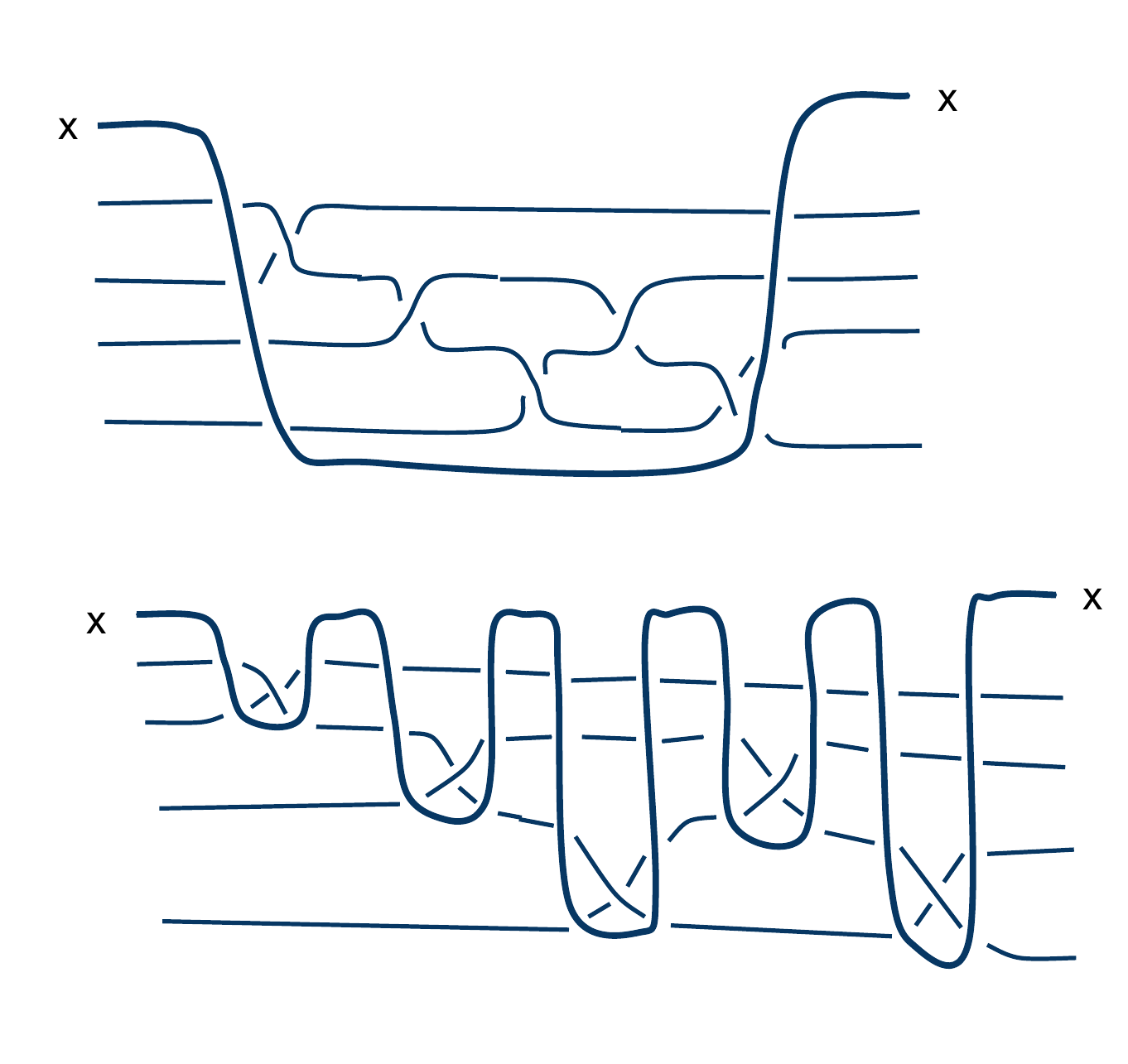}}

Therefore,  as we go to $0$ with the scale parameter $\varepsilon$  of the relative dilations with 
respect to $(x,\varepsilon)$, it is the same if we do it for the whole diagram encircled by a 
curve decorated with $(x,\varepsilon)$ or for the diagram where each crossing is encircled 
with $(x,\varepsilon)$. 

For encircled crossings the Reidemeister 3 move is valid. Being in an emergent algebra this means 
that we may use the Reidemeister 3 move in the interior of any encircled braid diagram, with the price 
of an error  $\mathcal{O}(\varepsilon)$ in the output of the encircled diagram. 
(Properly speaking, in order to have quantitative statements, we need to introduce a distance 
function on $X$, so we need to work with dilation structures \cite{buligadil1}.)

\centerline{\includegraphics[width=70mm]{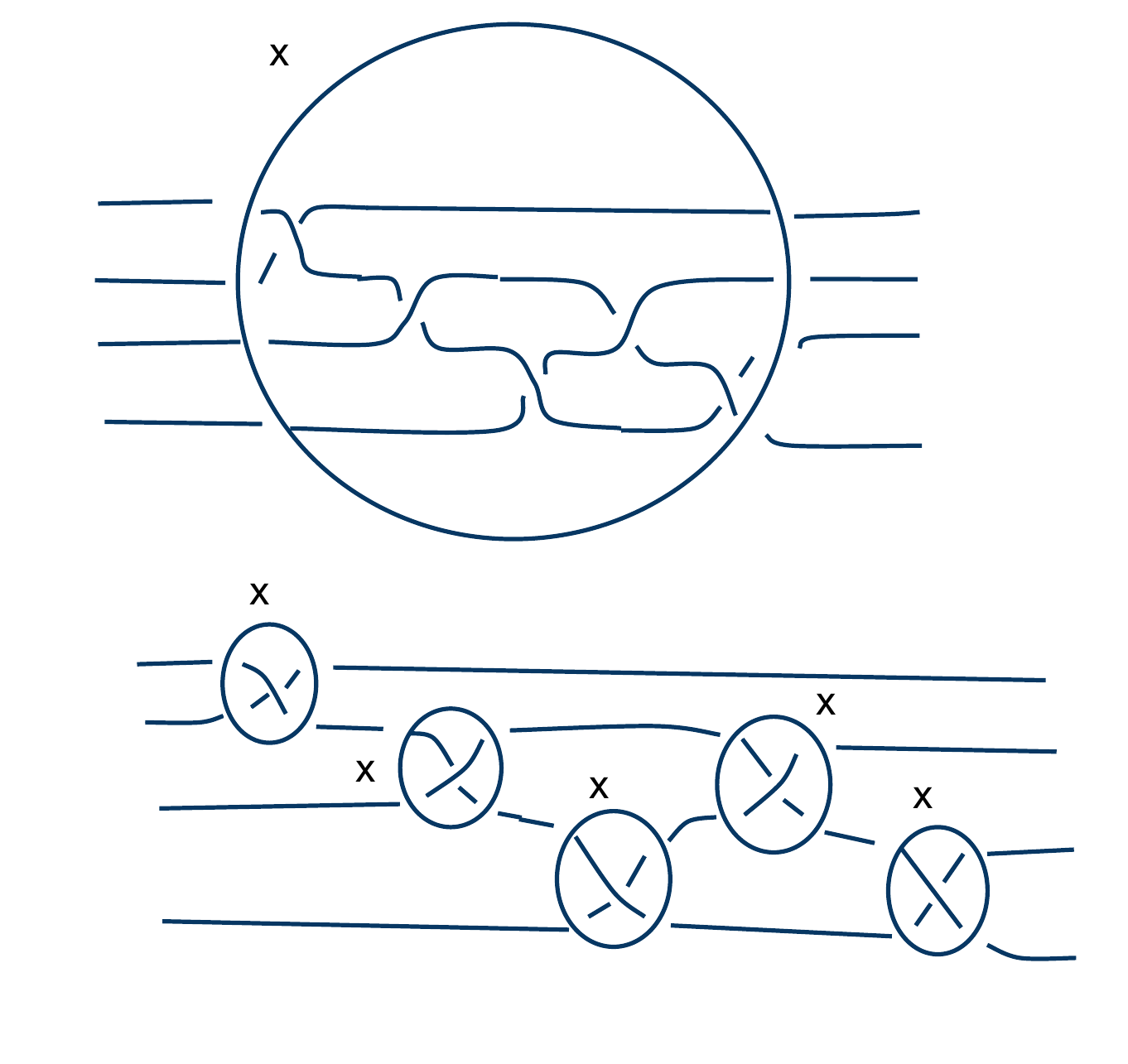}}

In the real world $\varepsilon$ does not go to $0$, it is only considered small, such that errors $\mathcal{O}(\varepsilon)$ are negligible.

What happens with the variable $f(u) \in \mathbb{R}$ sent through the channel $u \in X$, seen from the relative point of $x$? It is just differentiated, as shown in the next figure. 

\centerline{\includegraphics[width=120mm]{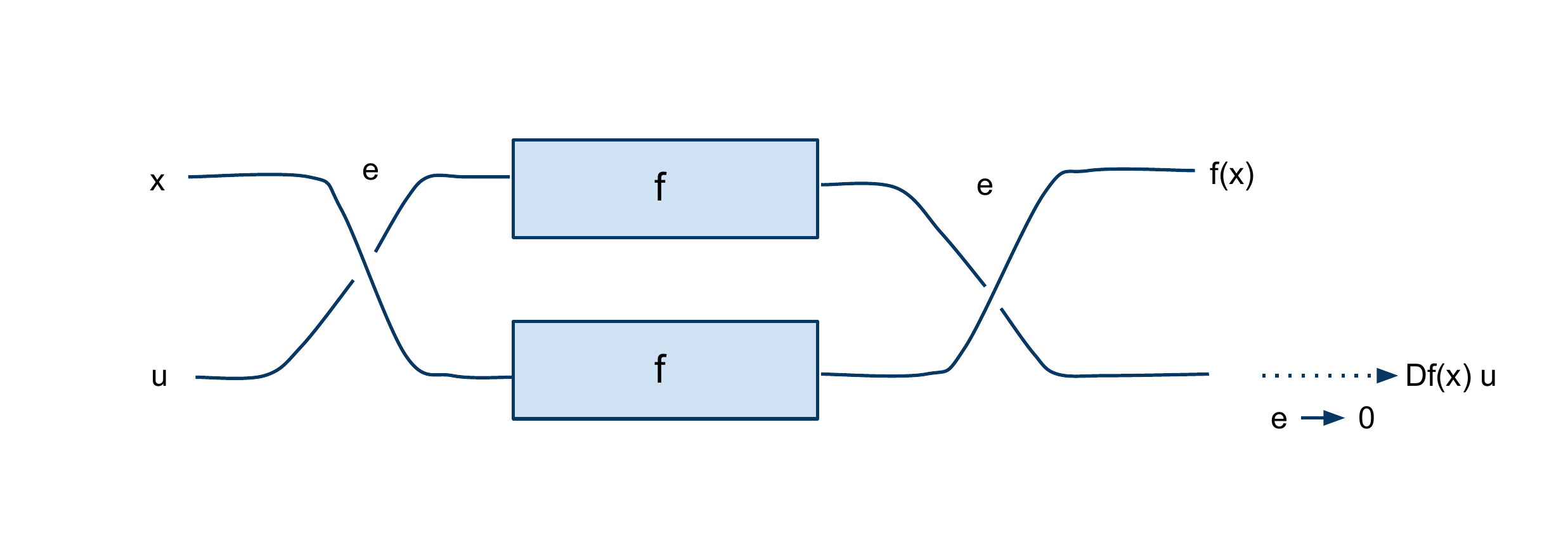}}

Here the two crossings in the diagram have different functions. The figure tells us that seen from the 
channel $x$, the first crossing from the left serves to define a new channel $\displaystyle x \circ_{\varepsilon} u$. Variables $f(x)$ and $\displaystyle f( x \circ_{\varepsilon} u)$ are communicated through their respective channels and we arrive at the second crossings, where we use dilation operation 
$\mathbb{R}$. As $\varepsilon$ goes to $0$, the output of this diagram is equivalent with the differentiation of $f$ in $x$, along the "direction" $u$. 

\section{Conclusion}

A decorated braids diagram represents an (exterior) process. This process is  seen from a visual point $x \in X$; the "seeing" is encoded as the functor of encircling diagrams. Visual processing 
is then described by differentiation of the stimuli (variables sent through channels). 

This is a very schematic and naive model, but one which I believe deserves further study.  All comments are welcome.

\paragraph{Acknowledgements.}  I  wish to 
express my thanks to Institut des Hautes Etudes Scientifiques, Paris, where a part of this work
has been done during a visit in  2010 and to Universidade Federal do Rio de Janeiro, where 
this work has been concluded.

\section{Appendix}

\subsection{Absolutes}

Let  $\Gamma$ be a topological commutative group. We suppose that 
$\Gamma$ as a topological space is separable.

\begin{definition}
Let $(X,\tau)$ be a topological space. $\tau$ is the collection of open sets in
$X$. A filter in $(X,\tau)$ is a function $\mu: \tau \rightarrow \left\{0,1\right\}$ such that: 
\begin{enumerate}
\item[(a)] $\mu(X) \ = \ 1$, 
\item[(b)] for any $A,B \in \tau$, if $A \subset B$ then $\mu(A) \leq \mu(B)$,
\item[(c)] for any $A,B \in \tau$ we have $\displaystyle \mu(A\cup B) + \mu(A\cap B) \geq \mu(A) + \mu(B)$.
\end{enumerate}
An absolute of  a separable topological commutative group $\Gamma$ is a
class $Abs(\Gamma)$  of filters $\mu$ in $\Gamma$ with the properties: 
\begin{enumerate}
\item[(i)] for any $\varepsilon \in \Gamma$ there are $A \in \tau$ and $\mu 
\in Abs(\Gamma)$ such that $\mu(A) = 1$ and $x \not \in A$, 
\item[(ii)] for any $\mu, \mu' \in Abs(\Gamma)$ there is $A \in \tau$ such that 
$\mu(A) > \mu'(A)$, 
\item[(iii)] for any $\varepsilon \in \Gamma$ and $\mu \in Abs(\Gamma)$ the 
transport of $\mu$ by $\varepsilon$, defined as $\varepsilon \, \mu (A) = 
\mu(\varepsilon A)$, belongs to $Abs(\Gamma)$. 
\end{enumerate}
Let $f: \Gamma \rightarrow (X,\tau)$ be a function from $\Gamma$ to a
separable topological space, let  $Abs(\Gamma)$ be an absolute of $\Gamma$, and 
$\mu \in Abs(\Gamma)$. We say that $f$ converges to $z \in X$ as 
$\varepsilon$ goes to $\mu$ if for any open set $A$ in $X$ with $z \in A$ we have
$\displaystyle \mu(f^{-1}(A)) = 1$. We write: 
$$\lim_{\varepsilon \rightarrow \mu} f(\varepsilon) \, = \, z$$
\label{dfilta}
\end{definition}

For example, if $\Gamma = (0,+\infty)$ with multiplication, then 
$\displaystyle Abs(\Gamma) = \left\{ 0 \right\}$ is an absolute, where "$0$" is the filter
defined by $0(A) = 1$ if and only if the number $0$ belongs to the closure of 
$A$ in $\mathbb{R}$. Also, $Abs(\Gamma) = \left\{ 0, \infty \right\}$ is an
absolute, where "$\infty$" is the filter defined by: $\infty(A) = 1$ if and only
if $A$ is unbounded. 

Let $\Gamma$ be a commutative separable topological group,  
 $\chi: \Gamma \rightarrow (0,+\infty)$ a continuous morphism and $Abs((0,+\infty))$ an 
 absolute of $(0,+\infty)$. Let  $Abs(\Gamma)$ be the class of filters on 
 $\Gamma$ constructed like this: $\mu \in Abs(\Gamma)$ if there exists 
 $\alpha \in Abs((0,+\infty))$ such that for any open set $A$ in $\Gamma$, 
 $\mu(A) = 1$ if there is an open set $B \subset (0,+ \infty)$  with 
 $\displaystyle \chi^{-1}(B) \subset A$ and $\alpha(B) = 1$. Then $Abs(\Gamma)$
 is an absolute of $\Gamma$. 
 
 Another example: let $\displaystyle \Gamma_{0}$ be a topological separable 
 commutative group,let $G$ be a finite commutative group and let 
 $\displaystyle \Gamma = \Gamma_{0} \times G$. We think now about $G$ and 
 $\displaystyle \Gamma_{0}$ as being 
 subgroups of $\Gamma$.  Let $\displaystyle
 Abs(\Gamma_{0})$ be an absolute of $\displaystyle \Gamma_{0}$. We construct 
 $Abs(\Gamma)$ as the collection of all filters $\mu$ on $\Gamma$ such that 
 there is $g \in G$ with $\displaystyle g \mu \in Abs(\Gamma_{0})$. Then
 $Abs(\Gamma)$ is an absolute of $\Gamma$.

\subsection{Conical groups}

\begin{definition}
A contractible group is a pair $(G,\alpha)$, where $G$ is a  
topological group with neutral element denoted by $e$, and $\alpha \in Aut(G)$ 
is an automorphism of $G$ such that: 
\begin{enumerate}
\item[-] $\alpha$ is continuous, with continuous inverse, 
\item[-] for any $x \in G$ we have the limit $\displaystyle 
\lim_{n \rightarrow \infty} \alpha^{n}(x) = e$. 
\end{enumerate}
\label{defunu}
\end{definition}

If $(G,\alpha)$ is a contractible group then $(G,\circ)$ is a irq, with: 
$$x \circ y \,  =  \, x \alpha(x^{-1} y)$$

Contractible groups are particular examples of conical groups. In
\cite{buligairq} we proved that conical groups, as well as some symmetric
spaces,  can be described as emergent algebras, coming from uniform idempotent
right quasigroups.

Conical groups are  particular cases of  group with dilations,  introduced in \cite{buliga2}, \cite{buligadil1}; we describe them further.

Let $\Gamma$ be a topological commutative groups with an absolute $Abs(\Gamma)$.

\begin{definition}
A group with dilations $(G,\delta)$ is a  topological group $G$  with  an action 
of $\Gamma$ (denoted by $\delta$), on $G$ such that for any $\mu \in Abs(\Gamma)$
\begin{enumerate}
\item[H0.] the limit  $\displaystyle \lim_{\varepsilon \rightarrow \mu} 
\delta_{\varepsilon} x  =  e$ exists and is uniform with respect to $x$ in a compact neighbourhood of the identity $e$.
\item[H1.] the limit
$$\beta (x,y)  =  \lim_{\varepsilon \rightarrow \mu} \delta_{\varepsilon}^{-1}
\left((\delta_{\varepsilon}x) (\delta_{\varepsilon}y ) \right)$$
is well defined in a compact neighbourhood of $e$ and the limit is uniform.
\item[H2.] the following relation holds
$$ \lim_{\varepsilon \rightarrow \mu} \delta_{\varepsilon}^{-1}
\left( ( \delta_{\varepsilon}x)^{-1}\right)  =  x^{-1}$$
where the limit from the left hand side exists in a neighbourhood of $e$ and is uniform with respect to $x$.
\end{enumerate}
\label{defgwd}
\end{definition}

\begin{definition}
A conical group $(N, \delta)$ is a   group with dilations  such that for any $\varepsilon \in \Gamma$  the dilation 
 $\delta_{\varepsilon}$ is a group morphism. 
\end{definition}

A conical group is the infinitesimal version of a group with 
dilations (\cite{buligadil1} proposition 2).

\begin{proposition}
Under the hypotheses H0, H1, H2,  $\displaystyle (G,\beta, \delta)$,  is a conical group, with operation 
$\displaystyle \beta$,  dilations $\delta$.
\label{here3.4}
\end{proposition}

One particular case is the one of contractible groups, definition \ref{defunu}, 
which are also normed groups. Indeed, in this case we may take 
$\Gamma = \mathbb{Z}$.

Locally compact conical groups are  locally compact groups admitting 
a contractive automorphism group. We begin with  the
definition of a contracting automorphism group \cite{siebert}, definition 5.1. 

\begin{definition}
Let $G$ be a locally compact group. An automorphism group on $G$ is a family 
$\displaystyle T= \left( \tau_{t}\right)_{t >0}$ in $Aut(G)$, such that 
$\displaystyle \tau_{t} \, \tau_{s} = \tau_{ts}$ for all $t,s > 0$. 

The contraction group of $T$ is defined by 
$$C(T) \ = \ \left\{ x \in G \mbox{ : } \lim_{t \rightarrow 0} \tau_{t}(x) = e
\right\} \quad .$$
The automorphism group $T$ is contractive if $C(T) = G$. 
\end{definition}

Next is proposition 5.4 \cite{siebert}, which gives a description of locally
compact groups which admit a contractive automorphism group. 

\begin{proposition}
For a locally compact group $G$ the following assertions are equivalent: 
\begin{enumerate}
\item[(i)] $G$ admits a contractive automorphism group;
\item[(ii)] $G$ is a simply connected Lie group whose Lie algebra admits a 
positive graduation.
\end{enumerate}
\label{psiebert}
\end{proposition}

The proof of the next proposition is an easy  application of 
the previously explained facts. 

\begin{proposition}
Let $(G,\delta)$ be a locally compact conical group. Then the associate 
irq  $(G,\circ)$ is an uniform irq. 
\label{pfirstop}
\end{proposition}

A particular class of locally compact groups which admit a contractive
automorphism group is  made by Carnot groups. They are related to sub-riemannian or 
Carnot-Carath\'eodory geometry, which is the study of non-holonomic manifolds
endowed with a Carnot-Carath\'eodory distance. Non-holonomic spaces were
discovered in 1926 by  G. Vr\u anceanu \cite{vra1},
\cite{vra2}. The  Carnot-Carath\'eodory distance on a non-holonomic space is 
inspired by Carath\'eodory \cite{cara} work from 1909   on the mathematical formulation of
thermodynamics. Such spaces appear in applications
to thermodynamics, to the mechanics of non-holonomic systems, in the study of
hypo-elliptic operators cf. H\"ormander \cite{hormander}, in harmonic analysis on
homogeneous cones cf. Folland, Stein \cite{fostein}, and
as boundaries of CR-manifolds.

The following result is a slight modification of  \cite{buligairq}, theorem 6.1,
consisting in the replacement of "contractible" by "conical" in the statement of
the theorem. 

\begin{theorem}
Let $(G,\alpha)$ be a locally compact conical  group and $G(\alpha)$ be the 
associated uniform irq. Then the irq is distributive, namely it satisfies the
relation: for any $\varepsilon, \lambda \in \Gamma$ 
\begin{equation}
x  \circ_{\varepsilon} \left( y \circ_{\lambda} z \right) \, = \, \left( x
\circ_{\varepsilon} y \right) \circ_{\lambda} \left( 
x \circ_{\varepsilon} z \right) 
\label{distributive}
\end{equation}

Conversely, let $(G, \circ)$ be a distributive uniform irq. Then there 
is a group operation on $G$ (denoted multiplicatively), with neutral element
$e$, such that:
\begin{enumerate}
\item[(i)] $\displaystyle xy \, = \, \Sigma^{e}(x, y)$ for any $x, y \in G$, 
\item[(ii)] for any $x, y \in G$ we have $\displaystyle 
x \circ_{\varepsilon} y \, = \, x (e \circ_{\varepsilon} (x^{-1} y))$. 
\end{enumerate}
In conclusion  there is a bijection between distributive $\Gamma$-uniform irqs and 
conical groups. 
\label{pgroudlin}
\end{theorem}

\subsection{Normed uniform irqs are dilation structures} 
\label{dilats}

For simplicity we shall list the   axioms of  a dilation structure $(X,d,\delta)$ without concerning about
domains and codomains of dilations. For the full definition of dilation
structures, as well as for their main properties and examples, see
\cite{buligadil1}, \cite{buligadil2}, \cite{buligasr}. The notion appeared from
my efforts to understand  the last section of the paper  
\cite{bell} (see also \cite{pansu}, \cite{gromovsr}, \cite{marmos1},
\cite{marmos2}).

Let $\Gamma$ be a topological commutative groups with an absolute $Abs(\Gamma)$ 
and with a morphism $\mid \cdot \mid : \Gamma \rightarrow (0,+\infty)$ such that
for any $\mu \in Abs(\Gamma)$ 
$$\lim_{\varepsilon \rightarrow \mu} \mid \varepsilon \mid \, = \, 0$$

\begin{definition}
A triple $(X, d, \delta)$ is a dilation structure if $(X, d)$ is a locally 
compact metric space and the dilation field 
 $$\delta: \Gamma \times \left\{ (x,y) \in X \times X \mbox{ : } y \in
 dom(\varepsilon, x) \right\} \rightarrow X \quad , \quad \delta(\varepsilon, x,
 y) \, = \, \delta^{x}_{\varepsilon} y $$
 gives to $X$ the structure of a uniform idempotent right quasigroup over 
 $\Gamma$, with the operation: for any $\varepsilon \in \Gamma$  
 $$x \, \circ_{\varepsilon} \, y \, = \, \delta_{\varepsilon}^{x} y$$
 Moreover, the distance is compatible with the dilations, in the sense: 
 
\begin{enumerate}
\item[A1.]  the uniformity on $(X,\delta)$ is the one induced by the
distance $d$,

\item[A2.] There is $A > 1$ such that  for any 
$x$ there exists
 a  function $\displaystyle (u,v) \mapsto d^{x}(u,v)$, defined for any
$u,v$ in the closed ball (in distance d) $\displaystyle
\bar{B}(x,A)$, such that for any $\mu \in Abs(\Gamma)$
$$\lim_{\varepsilon \rightarrow \mu} \quad \sup  \left\{  \mid \frac{1}{\mid 
\varepsilon \mid} d(\delta^{x}_{\varepsilon} u,
\delta^{x}_{\varepsilon} v) \ - \ d^{x}(u,v) \mid \mbox{ :  } u,v \in \bar{B}_{d}(x,A)\right\} \ =  \ 0$$
uniformly with respect to $x$ in compact set. 
Moreover the uniformity induced by $d^{x}$ is the same as the uniformity induced
by $d$, in particular $\displaystyle d^{x}(u,v) = 0$ implies $u = v$.

\end{enumerate}

\label{defweakstrong}
\end{definition}

The conclusion is therefore that adding a distance in the story of uniform 
irqs gives us the notion of a dilation structure. 

We go a bit into details. 

\begin{proposition}
Let $(X,d, \delta)$ be a dilation structure,  $x \in X$, and let 
$$\delta^{x}_{\varepsilon} \, d (u, v) \, = \, \frac{1}{\mid \varepsilon \mid}
\, d( \delta^{x}_{\varepsilon} u , \delta^{x}_{\varepsilon} v )$$
 Then the net of metric spaces $\displaystyle (\bar{B}_{d}(x,A),
 \delta^{x}_{\varepsilon} d)$ converges in the Gromov-Hausdorff sense to the 
 metric space $\displaystyle (\bar{B}_{d}(x,A), d^{x})$. Moreover this metric
 space is a metric cone, in the following sense: for any $\lambda \in \Gamma$ 
  we have
 $$d^{x} ( \delta^{x}_{\lambda} u , \delta^{x}_{\lambda} v ) \, = \, 
 \mid \lambda \mid 
 \, d^{x}(u,v)$$
 \label{pizo}
 \end{proposition}
 
 \paragraph{Proof.} 
 The first part of the proposition is just a reformulation  of axiom A2, 
 without the condition of uniform convergence. For the second part 
 remark  that 
  
 $$\frac{1}{\mid \varepsilon \mid} d( \delta^{x}_{\varepsilon} \, 
 \delta^{x}_{\lambda} u , \delta^{x}_{\varepsilon} \, 
 \delta^{x}_{\lambda} v ) \, = \, \mid \lambda \mid \, \frac{1}{\mid 
 \varepsilon \lambda \mid} d( \delta^{x}_{\varepsilon \lambda}  u , 
 \delta^{x}_{\varepsilon \lambda}  v ) $$
 Therefore if we pass to the limit 
 with $\varepsilon \rightarrow \mu$ in these two relations we get the desired
 conclusion. \quad $\square$

Particular examples of dilation structures are given by normed groups with
dilations. 

\begin{definition} A normed group with dilations $(G, \delta, \| \cdot \|)$ is a 
group with dilations  $(G, \delta)$ endowed with a continuous norm  
function $\displaystyle \|\cdot \| : G \rightarrow \mathbb{R}$ which satisfies 
(locally, in a neighbourhood  of the neutral element $e$) the properties: 
 \begin{enumerate}
 \item[(a)] for any $x$ we have $\| x\| \geq 0$; if $\| x\| = 0$ then $x=e$, 
 \item[(b)] for any $x,y$ we have $\|xy\| \leq \|x\| + \|y\|$, 
 \item[(c)] for any $x$ we have $\displaystyle \| x^{-1}\| = \|x\|$, 
 \item[(d)] the limit 
$\displaystyle \lim_{\varepsilon \rightarrow \mu} \frac{1}{\mid\varepsilon \mid} \| \delta_{\varepsilon} x \| = \| x\|^{N}$ 
 exists, is uniform with respect to $x$ in compact set, 
 \item[(e)] if $\displaystyle \| x\|^{N} = 0$ then $x=e$.
  \end{enumerate}
  \label{dnco}
  \end{definition}

In a normed group with dilations we have a natural left invariant distance given by
\begin{equation}
d(x,y) = \| x^{-1}y\| \quad . 
\label{dnormed}
\end{equation}
Any normed group with dilations has an associated dilation structure on it.  In a group with dilations $(G, \delta)$  we define dilations based in any point $x \in G$ by 
 \begin{equation}
 \delta^{x}_{\varepsilon} u = x \delta_{\varepsilon} ( x^{-1}u)  . 
 \label{dilat}
 \end{equation}

The following result is theorem 15 \cite{buligadil1}. 

\begin{theorem}
Let $(G, \delta, \| \cdot \|)$ be  a locally compact  normed  group with dilations. Then $(G, d, \delta)$ is 
a dilation structure, where $\delta$ are the dilations defined by (\ref{dilat}) and the distance $d$ is induced by the norm as in (\ref{dnormed}). 
\label{tgrd}
\end{theorem}

The general theorem \ref{mainthm} has a stronger conclusion in the case of
dilation structures, namely "conical groups" are replaced by "normed conical
groups".

\end{document}